\documentclass[aap,MSNbibl,seceqn,dvips]{arximspdf}
\usepackage{url,breakurl}
\doi{10.1214/13-AAP963} 
\volume{24}
\issue{5}
\pubyear{2014}
\firstpage{1803}
\lastpage{1834}

\makeatletter
\def\cal{\mathcal}
\newcommand{\eqref}[1]{(\ref{#1})}

\newtheorem{Theorem}{Theorem}[section]
\newproclaim{Definition}[Theorem]{Definition}
\newtheorem{Proposition}[Theorem]{Proposition}

\newproclaim{Assumption}[Theorem]{Assumption}

\newtheorem{Lemma}[Theorem]{Lemma}

\newproclaim{Remark}[Theorem]{Remark}
\newproclaim{Example}[Theorem]{Example}

\def\E{\mathbb{E}}
\def\F{\mathbb{F}}
\def\P{\mathbb{P}}
\def\R{\mathbf{R}}
\def\R{\mathbb{R}}
\def\N{\mathbb{N}}

\def\Ac{{\cal A}}
\def\Ec{{\cal E}}
\def\Fc{{\cal F}}
\def\Lc{{\cal L}}
\def\Mc{{\cal M}}
\def\Pc{{\cal P}}
\def\Sc{{\cal S}}
\def\Uc{{\cal U}}

\def\Lb{{\bar L}}
\def\Pb{\overline{\P}}
\def\Xb{\overline{X}}
\def\Ub{\overline{U}}
\def\Ab{\overline{A}}
\def\Bb{\overline{B}}
\def\Gb{\overline{G}}

\def\Mbf{\mathbf{M}}

\def\Om{\Omega}
\def\Omb{\overline{\Omega}}
\def\eps{\varepsilon}

\def\Fcb{\overline{{\cal F}}}
\def\Fbb{\overline{\mathbb{F}}}

\def\Pcb{\overline{\Pc}}

\def\nub{\overline{\nu}}
\def\mub{\overline{\mu}}
\def\sigmab{\overline{\sigma}}
\def\Lb{\overline{L}}
\def\Phib{\overline{\Phi}}

\makeatother

\begin{document}
\begin{frontmatter}

\title{Discrete-time probabilistic approximation of
path-dependent stochastic control problems}
\runtitle{Monte Carlo for stochastic control}

\begin{aug}
\author{\fnms{Xiaolu} \snm{Tan}\corref{}\ead[label=e1]{tan@ceremade.dauphine.fr}}
\runauthor{X. Tan}
\affiliation{Ceremade, University of Paris-Dauphine}
\address{Ceremade, Universite Paris-Dauphine\\
Place du Marechal de Lattre de Tassigny\\
75775, Paris Cedex 16\\
France\\
\printead{e1}}
\end{aug}

\received{\smonth{8} \syear{2012}}
\revised{\smonth{8} \syear{2013}}

%
\begin{abstract}
We give a probabilistic interpretation of the Monte Carlo scheme
proposed by Fahim, Touzi and Warin [\textit{Ann. Appl. Probab.} \textbf{21}
(2011) 1322--1364] for fully nonlinear parabolic PDEs, and hence
generalize it to the path-dependent (or non-Markovian) case for a
general stochastic control problem. A general convergence result is
obtained by a weak convergence method in the spirit of Kushner and
Dupuis
[\textit{Numerical Methods for Stochastic Control Problems in Continuous Time}
(1992) Springer]. We also get a rate of
convergence using the invariance principle technique as in Dolinsky
[\textit{Electron. J. Probab.} \textbf{17} (2012) 1--5], which is better
than that obtained by viscosity solution method. Finally, by
approximating the conditional expectations arising in the numerical
scheme with simulation-regression method, we obtain an implementable scheme.
\end{abstract}

%

\end{frontmatter}

\section{Introduction}
Stochastic optimal control theory is largely applied in economics,
finance, physics and management problems. Since its development,
numerical methods for stochastic control problems have also been
largely investigated. For the Markovian control problem, the value
function can usually be characterized by Hamilton--Jacob--Bellman
(HJB)
equations, then many numerical methods are also given as numerical
schemes for PDEs. In this context, a powerful tool to prove the
convergence is the monotone convergence of viscosity solution method of
Barles and Souganidis \cite{BarlesSouganidis}.

In the one-dimensional case, the explicit finite difference scheme can
be easily constructed and implemented, and the monotonicity is
generally guaranteed under the Courant--Friedrichs--Lewy (CFL)
condition. In two dimensional cases, Bonnans, Ottenwaelter and Zidani
\cite{Bonnans2004} proposed a numerical algorithm to construct
monotone explicit schemes. Debrabant and Jakobsen \cite
{DebrabantJakobsen} gave a semi-Lagrangian scheme which is easily
constructed to be monotone but needs finite difference grid together
with interpolation method for the implementation. In general, these
methods may be relatively efficient in low dimensional cases; while in
high dimensional cases, a Monte Carlo method is preferred if possible.

As a generalization of the Feynman--Kac formula, the backward
stochastic differential equation (BSDE) opens a way for the Monte Carlo
method for optimal control problems; see, for example, Bouchard and
Touzi \cite{Bouchard2004}, Zhang \cite{Zhang2004}. Generally speaking,
the BSDE covers the controlled diffusion processes problems of which
only the drift part is controlled. However, it cannot include the
general control problems when the volatility part is also controlled.
This is one of the main motivations of recent developments of second
order BSDE (2BSDE) by Cheridito, Soner, Touzi and Victoir \cite
{Cheridito2007} and Soner, Touzi and Zhang \cite{SonerTouziZhang1}.
Motivated by the 2BSDE theory in \cite{Cheridito2007}, and also
inspired by the numerical scheme of BSDEs, Fahim, Touzi and Warin \cite
{FahimTouziWarin} proposed a probabilistic numerical scheme for fully
nonlinear parabolic PDEs. In their scheme, one needs to simulate a
diffusion process and estimate the value function as well as the
derivatives of the value function arising in the PDE by conditional
expectations, and then compute the value function in a backward way on
the discrete time grid. The efficiency of this Monte Carlo scheme has
been shown by several numerical examples, and we refer to Fahim, Touzi
and Warin \cite{FahimTouziWarin}, Guyon and Henry-Labord\`ere \cite
{GuyonHenryLabordere} and Tan \cite{TanSplitting} for the implemented
examples.

However, instead of probabilistic arguments, the convergence of this
scheme is proved by techniques of monotone convergence of viscosity
solution of Barles and Souganidis \cite{BarlesSouganidis}. Moreover,
their scheme can be only applied in the Markovian case when the value
function is characterized by PDEs.

The main contribution of this paper is to give a probabilistic
interpretation to the Monte Carlo scheme of Fahim, Touzi and Warin
\cite
{FahimTouziWarin} for fully nonlinear PDEs, which allows one to
generalize it to the non-Markovian case for a general stochastic
optimal control problem. One of the motivations for the non-Markovian
generalization comes from finance to price the path-dependent exotic
derivative options in the uncertain volatility model.

Our general convergence result is obtained by weak convergence
techniques in spirit of Kushner and Dupuis \cite{KushnerDupuis}. In
contrast to \cite{KushnerDupuis}, where the authors define their
controlled Markov chain in a descriptive way, we give our controlled
discrete-time semimartingale in an explicit way using the cumulative
distribution functions. Moreover, we introduce a canonical space for
the control problem following El Karoui, H{\.u}{\.u} Nguyen and Jeanblanc \cite{ElKaroui1987}, which allows us to explore the convergence conditions
on the reward functions. We also provide a convergence rate using the
invariance principle techniques of Sakhanenko \cite{Sakhanenko} and
Dolinsky~\cite{Dolinsky}. Compared to the scheme in \cite{Dolinsky} in
the context of $G$-expectation (see, e.g., Peng~\cite{PengGexpec} for
$G$-expectation), our scheme is implementable using
simulation-regression method.

The rest of the paper is organized as follows. In Section~\ref{sec:numscheme}, we first introduce a general path-dependent stochastic
control problem and propose a numerical scheme. Then we give the
assumptions on the diffusion coefficients and the reward functions, as
well as the main convergence results, including the general convergence
and a rate of convergence. Next in Section~\ref{sec:proba_interp}, we
provide a probabilistic interpretation of the numerical scheme, by
showing that the numerical solution is equivalent to the value function
of a controlled discrete-time semimartingale problem. Then we complete
the proofs of the convergence results in Section~\ref{sec:conv_proof}.
Finally, in Section~\ref{sec:implem_scheme}, we discuss some issues
about the implementation of our numerical scheme, including a
simulation-regression method.

\textit{Notation}. We denote by $S_d$ the space of all $d
\times d$
matrices, and by $S_d^+$ the space of all positive symmetric $d \times
d$ matrices.
Given a vector or a matrix $A$, then $A^{\top}$ denotes its transposition.
Given two $d \times d$ matrix $A$ and $B$, their product is defined by
$A \cdot B:= \operatorname{Tr}(A B^{\top})$ and $|A|:= \sqrt{A \cdot A}$.
Let $\Om^d:= C([0,T], \R^d)$ be the space of all continuous paths
between $0$ and $T$, denote $|\mathbf{x}|:= \sup_{0 \le t \le T}
|\mathbf{x}_t|$ for
every $\mathbf{x}\in\Om^d$.
In the paper, $E$ is a fixed compact Polish space, we denote
\[
Q_T:=[0,T] \times\Om^d \times E.
\]
Suppose that $(X_{t_k})_{0 \le k \le n}$ is a process defined on the
discrete time grid $(t_k)_{0 \le k \le n}$ of $[0,T]$ with $t_k:= k h$
and $h:= \frac{T}{n}$. We usually write it as $(X_k)_{0 \le k \le n}$,
and denote by $\widehat{X}$ its linear interpolation path on $[0,T]$.
In the paper, $C$ is a constant whose value may vary from line to line.

\section{A numerical scheme for stochastic control problems}
\label{sec:numscheme}

\subsection{A path-dependent stochastic control problem}
\label{subsec:ContrPb}

Let $(\Om, \Fc, \P)$ be a complete probability space containing a
$d$-dimensional standard Brownian motion $W$, $\F=(\Fc_t)_{0 \le t
\le
T}$ be the natural Brownian filtration. Denote by $\Om^d:= C([0,T],
\R
^d)$ the space of all continuous paths between $0$ and $T$. Suppose
that $E$ is a compact Polish space with a complete metric $d_E$,
$(\sigma, \mu)$ are bounded continuous functions defined on $Q_T:=
[0,T] \times\Om^d \times E$ taking value in $S_d \times\R^d $.
We fix a constant
$x_0 \in\R^d$ through out the paper. Then given a $\F$-progressively
measurable $E$-valued process $\nu= (\nu_t)_{0 \le t \le T}$, denote
by $X^{\nu}$ the controlled diffusion process which is the strong
solution to
%
\begin{equation}
\label{eq:SDE} X^{\nu}_t = x_0 + \int
_0^t \mu\bigl(s, X^{\nu}_{\cdot},
\nu_s\bigr) \,ds + \int_0^t \sigma
\bigl(s, X^{\nu}_{\cdot}, \nu_s\bigr)
\,dW_s.
\end{equation}
To ensure the existence and uniqueness of the strong solution to the
above equation \eqref{eq:SDE}, we suppose that for every progressively
measurable process $(X, \nu)$, the processes $\mu(t,X_{\cdot}, \nu_t)$
and $\sigma(t,X_{\cdot},\nu_t)$ are progressively measurable. In
particular, $\mu$ and $\sigma$ depend on the past trajectory of $X$.
Further, we suppose that there is some constant $C$ and a continuity
module $\rho$, which is an increasing function on $\R^+$ satisfying
$\rho(0^+) = 0$, such that
%
\begin{eqnarray}
\label{eq:cond_mu_sigma} && \bigl| \mu(t_1, \mathbf{x}_1,
u_1) - \mu(t_2, \mathbf{x}_2,
u_2) \bigr| + \bigl| \sigma(t_1, \mathbf{x}_1,
u_1) - \sigma(t_2, \mathbf{x}_2,
u_2) \bigr|
\nonumber
\\[-8pt]
\\[-8pt]
\nonumber
&&\qquad\le C \bigl|\mathbf{x}_1^{t_1} - \mathbf{x}_2^{t_2}\bigr|
+ \rho\bigl(|t_1 - t_2|+d_E(u_1,
u_2)\bigr),
\end{eqnarray}
where for every $(t,\mathbf{x}) \in[0,T] \times\Om^d$, we denote
$\mathbf{x}^t_s:= \mathbf{x}
_s 1_{[0,t]}(s) + \mathbf{x}_t 1_{(t,T]}(s)$.
Let $\Phi\dvtx  \mathbf{x}\in\Om^d \to\R$ and $L \dvtx (t,\mathbf{x},u)
\in Q_T \to\R$
be the continuous reward functions, and denote by $\Uc$ the collection
of all $E$-valued $\F$-progressively measurable processes, the main
purpose of this paper is to approximate numerically the following
optimization problem:
%
\begin{equation}
\label{eq:V} V:= \sup_{\nu\in\Uc} \E \biggl[ \int
_0^T L \bigl(t,X^{\nu
}_{\cdot
},
\nu_t \bigr) \,dt + \Phi\bigl(X_{\cdot}^{\nu}\bigr)
\biggr].
\end{equation}
Similarly to $\mu$ and $\sigma$, we suppose that for every
progressively measurable process $(X, \nu)$, the process $t \mapsto
L(t, X_{\cdot}, \nu_t)$ is progressively measurable. Moreover, to
ensure that the expectation in \eqref{eq:V} is well defined, we shall
assume later that $L$ and $\Phi$ are of exponential growth in $\mathbf
{x}$ and
discuss their integrability in Proposition \ref{prop:Y_integ}.

\subsection{The numerical scheme}

In preparation of the numerical scheme, we shall fix, through out the
paper, a progressively measurable function $\sigma_0 \dvtx [0,T] \times
\Om^d
\to S_d$ such that
%
\begin{eqnarray}
\bigl| \sigma_0(t_1, \mathbf{x}_1) -
\sigma_0(t_2, \mathbf{x}_2)\bigr | \le C \bigl|
\mathbf{x} _1^{t_1} - \mathbf{x}_2^{t_2}\bigr|
+ \rho\bigl(|t_1 - t_2|\bigr)
\nonumber\\
 \eqntext{\forall(t_1, \mathbf{x}_1), (t_2,
\mathbf{x}_2) \in [0,T] \times\Om^d,}
\end{eqnarray}
and with $\eps_0 > 0$, $\sigma_0 \sigma_0^{\top}(t,\mathbf{x}) \ge
\eps_0
I_d$ for every $(t,\mathbf{x}) \in[0,T] \times\Om^d$.
Denote
%
\begin{eqnarray}
\label{eq:note_a_u} \sigma_0^{t,\mathbf{x}}&:=& \sigma_0(t,
\mathbf{x}),\qquad a_0^{t,\mathbf{x}}:= \sigma _0^{t,\mathbf{x}}
\bigl(\sigma_0^{t,\mathbf{x}}\bigr)^{\top},
\nonumber
\\[-8pt]
\\[-8pt]
\nonumber
a^{t,\mathbf{x}}_u&:=& \sigma\sigma^{\top}(t, \mathbf{x}, u)
- a^{t,\mathbf{x}}_0, \qquad b^{t,\mathbf{x}}_u:= \mu(t,
\mathbf{x},u).
\end{eqnarray}
Then we define a function $G$ on $[0,T] \times\Om^d \times S_d
\times\R^d$ by
%
\begin{equation}
\label{eq:functionG} G(t, \mathbf{x}, \gamma, p):= \sup_{u \in E}
\biggl( L(t,\mathbf {x}, u) + \frac
{1}{2} a^{t,\mathbf{x}}_u
\cdot\gamma+ b^{t,\mathbf{x}}_u \cdot p \biggr),
\end{equation}
which is clearly convex in $(\gamma, p)$ as the supremum of a family
of linear functions, and lower-semicontinuous in $(t,\mathbf{x})$ as the
supremum of a family of continuous functions. Let $n \in\N$ denote the
time discretization by $h:= \frac{T}{n}$ and $t_k:= h k$.

Let us take the standard $d$-dimensional Brownian motion $W$ in the
complete probability space $(\Om, \Fc, \P)$. For simplicity, we denote
$W_k:= W_{t_k}$,\break  $\Delta W_k:= W_k - W_{k-1}$, $\Fc^W_k:= \sigma
(W_0, W_1, \ldots, W_k)$ and $\E^W_k[\cdot]:= \E[ \cdot| \Fc^W_k]$.
Then we have a process $X^0$ on the discrete grid $(t_k)_{0 \le k \le
n}$ defined by
%
\begin{equation}
\label{eq:defX0} X^0_0:= x_0,\qquad
X^0_{k+1}:= X^0_k +
\sigma_0\bigl(t_k, \widehat {X}{}^0_{\cdot}
\bigr) \Delta W_{k+1},
\end{equation}
where $\widehat{X}{}^0$ denotes the linear interpolation process of
$(X^0_k)_{0 \le k \le n}$ on interval $[0,T]$.

Then, for every time discretization $h$, our numerical scheme is given by
%
\begin{equation}
\label{eq:NumSchem} Y^h_k:= \E^W_k
\bigl[Y^h_{k+1}\bigr] + h G\bigl(t_k,
\widehat{X}{}^0_{\cdot}, \Gamma^h_k,
Z^h_k \bigr),
\end{equation}
with terminal condition
%
\begin{equation}
\label{eq:NumTermCond} Y^h_n:= \Phi\bigl(\widehat{X}{}^0_{\cdot}
\bigr),
\end{equation}
where $G$ is defined by \eqref{eq:functionG} and
\begin{eqnarray*}
\Gamma^h_k&:=& \E^W_k \biggl[
Y^h_{k+1} \bigl(\sigma_{0,k}^{\top}
\bigr)^{-1} \frac
{\Delta W_{k+1} \Delta W^{\top}_{k+1} - h I_d}{h^2} \sigma_{0,k}^{-1}
\biggr],\\
 Z^h_k&:=& \E^W_k \biggl[
Y^h_{k+1} \bigl(\sigma_{0,k}^{\top}
\bigr)^{-1} \frac{\Delta W_{k+1}}{h} \biggr],
\end{eqnarray*}
with $\sigma_{0,k}:= \sigma_0^{t_k, \widehat{X}{}^0} = \sigma_0(t_k,
\widehat{X}{}^0_{\cdot})$.

\begin{Remark}
By its definition, $Y^h_k$ is a measurable function of $(X^0_0, \ldots
,\break  X^0_k)$. We shall show later in Proposition \ref{prop:Y_integ} that
the function $Y^h_k(x_0, \ldots, x_k)$ is of exponential growth in
$\max_{0 \le i \le k}|x_i|$ under appropriate conditions, and hence the
conditional expectations in \eqref{eq:NumSchem} are well defined.
Therefore, the above scheme \eqref{eq:NumSchem} should be well defined.
\end{Remark}

\begin{Remark} \label{rem:MarkovScheme}
In the Markovian case, when the function $\Phi(\cdot)$ [resp.,
$L(t,\cdot,u)$, $\mu(t,\cdot,u)$ and $\sigma(t,\cdot,u)$] only depends
on $X_T$ (resp., $X_t$), so that the function $G(t,\cdot,\gamma,z)$ only
depends on $X_t$ and the value function of the optimization problem
\eqref{eq:V} can be characterized as the viscosity solution of a
nonlinear PDE
\[
- \partial_t v - \tfrac{1}{2} a_0(t,x) \cdot
D^2 v - G\bigl(t,x, D^2 v, D v\bigr) = 0.
\]
Then the above scheme reduces to that proposed by Fahim, Touzi and
Warin~\cite{FahimTouziWarin}.
\end{Remark}

\subsection{The convergence results of the scheme}

Our main idea to prove the convergence of the scheme \eqref
{eq:NumSchem}, \eqref{eq:NumTermCond} is to interpret it as an
optimization problem on a system of controlled discrete-time\vadjust{\goodbreak}
semimartingales, which converge weakly to the controlled diffusion
processes. Therefore, a reasonable assumption is that $\Phi$ and $L$
are bounded continuous on $\Om^d$ [i.e., $\Phi(\cdot), L(t,\cdot,u)
\in
C_b(\Om^d)$], or they belong to the completion space of $C_b(\Om^d)$
under an appropriate norm. We shall suppose that $\Phi$ (resp., $L$) is
continuous in $\mathbf{x}$ [resp., $(t,\mathbf{x}, u)$], and there are
a constant $C$
and continuity modules $\rho_0$, $(\rho_N)_{N \ge1}$ such that for
every $(t,\mathbf{x},u) \in Q_T$ and $(t_1,\mathbf{x}_1)$,
$(t_2,\mathbf{x}_2) \in[0,T]\times
\Om^d$ and $N \ge1$,
%
\begin{equation}
\label{eq:condPhi1} \qquad
\cases{\bigl |\Phi( \mathbf{x})\bigr| +
\bigl|L(t,\mathbf{x},u)\bigr| \le C \exp \bigl(
C |\mathbf{x}| \bigr), \vspace*{2pt}
\cr
\bigl|L(t_1,\mathbf{x},u) -
L(t_2,\mathbf{x},u)\bigr| \le \rho_0\bigl(|t_1 -
t_2|\bigr),\vspace*{2pt}
\cr
\bigl|\Phi_N(\mathbf{x}_1)
- \Phi_N(\mathbf{x}_2)\bigl| + \bigl|L_N(t,\mathbf
{x}_1,u) - L_N(t,\mathbf{x}_2,u)\bigr| \le
\rho_N\bigl(|\mathbf{x}_1 - \mathbf{x}_2|\bigr), }
\end{equation}
where $\Phi_N:= (-N) \vee(\Phi\land N)$ and $ L_N:= (-N) \vee(L
\land N)$.

Denote
%
\begin{equation}
\label{eq:defmG} m_G:= \min_{(t,\mathbf{x}, u) \in Q_T,  w \in\R^d} \biggl(
\frac
{1}{2}w^{\top} a^{t,\mathbf{x}}_u w +
b^{t,\mathbf{x}}_u \cdot w \biggr)
\end{equation}
and
%
\begin{equation}
\label{eq:defh0} h_0:= 1_{m_G = 0} T + 1_{m_G < 0} \min
_{(t,\mathbf{x}, u) \in
Q_T} - m_G^{-1} \biggl(1 -
\frac{1}{2} a^{t,\mathbf{x}}_u \cdot \bigl(
a_0^{t,\mathbf{x}} \bigr)^{-1} \biggr),
\end{equation}
where $a^{t,\mathbf{x}}_0$, $a^{t,\mathbf{x}}_u$ and $b^{t,\mathbf
{x}}_u$ are defined in
\eqref{eq:note_a_u}. Clearly, $m_G \le0$.

\begin{Assumption} \label{assum:G}
For every $(t,\mathbf{x},u) \in Q_T$, we have
\[
a_u^{t,\mathbf{x}} \ge 0 \quad\mbox{and} \quad1 - \tfrac{1}{2}
a^{t,\mathbf{x}}_u \cdot \bigl( a_0^{t,\mathbf{x}}
\bigr)^{-1} \ge 0.
\]
Further, the constants $m_G > - \infty$ and $h_0 > 0$.
\end{Assumption}

\begin{Remark} \label{rem:positive_density}
Assumption \ref{assum:G} is almost equivalent to Assumption F of \cite
{FahimTouziWarin} in the context of the control problem, and it
implies that the drift $\mu$ and $\sigma$ are uniformly bounded, as
assumed at the beginning of Section~\ref{subsec:ContrPb}. In
particular, it follows that when $m_G < 0$, we have
\[
1 - \tfrac{1}{2} a^{t,\mathbf{x}}_u \cdot \bigl(
a_0^{t,\mathbf{x}} \bigr)^{-1} + h m_G \ge 0\qquad
\mbox{for every } (t,\mathbf{x},u) \in Q_T \mbox { and } h \le
h_0.
\]
Moreover, since $a_0$ is supposed to be nondegenerate, the assumption
implies that $\sigma\sigma^{\top}(t,\mathbf{x},u)$ is nondegenerate
for all
$(t,\mathbf{x},u) \in Q_T$.
The nondegeneracy condition may be inconvenient in practice (see,
e.g., Example \ref{exam:UVM}), we shall also provide more discussions
and examples in Section~\ref{subsec:implem_G_expec}.
\end{Remark}

\begin{Remark}
When $a_u^{t,\mathbf{x}} \ge\eps I_d$ uniformly for some $\eps> 0$,
we get
immediately $m_G > - \infty$ since $b^{t,\mathbf{x}}_u$ is uniformly bounded.
When $a_u^{t,\mathbf{x}}$ degenerates, $m_G > -\infty$ implies that
$b_u^{t,\mathbf{x}
}$ lies in the image of $a_u^{t,\mathbf{x}}$.
\end{Remark}

\begin{Proposition} \label{prop:Y_integ}
Suppose that the reward functions $L$ and $\Phi$ satisfy \eqref
{eq:condPhi1}, then the optimal value $V$ in \eqref{eq:V} is finite.\vspace*{1pt}
Suppose in addition that Assumption \ref{assum:G} holds true. Then for
every fixed $n \in\N$ ($h:= \frac{T}{n}$) and every\vadjust{\goodbreak} $0 \le k \le n$,
as a function of $(X_0, \ldots, X_k)$, $Y_k^h(x_0, \ldots, x_k)$ is
also of exponential growth in $\max_{0 \le i \le k} |x_i|$. And hence
$Y^h_k$ is integrable in \eqref{eq:NumSchem}, the numerical scheme
\eqref{eq:NumSchem} is well defined.
\end{Proposition}
The proof is postponed until Section~\ref{subsec:technique_lemma} after a
technical lemma.

Our main results of the paper are the following two convergence
theorems, whose proofs are left in Section~\ref{sec:conv_proof}.

\begin{Theorem} \label{theo:convergence}
Suppose that $L$ and $\Phi$ satisfy \eqref{eq:condPhi1} and Assumption
\ref{assum:G} holds true. Then
\[
Y^h_0 \to V \qquad\mbox{as } h \to0.
\]
\end{Theorem}

To derive a convergence rate, we suppose further that $E$ is a compact
convex subset of $S^+_d \times\R^d$, and for every $(t,\mathbf{x},
u) =
(t,\mathbf{x},
a,b) \in Q_T$,
%
\begin{eqnarray}
\label{eq:condUplus} a > 0, \qquad\mu(t,\mathbf{x},u)& =& \mu(t,\mathbf{x},a,b) = b,
\nonumber
\\[-8pt]
\\[-8pt]
\nonumber
\sigma(t, \mathbf{x},u) &=& \sigma(t, \mathbf{x},a,b) = a^{1/2}.
\end{eqnarray}
Moreover, we suppose that
$L(t,\mathbf{x},u) = \ell(t,\mathbf{x}) \cdot u$ for some continuous
function $\ell\dvtx
[0,T] \times\Om^d \to S_d \times\R^d$ and that there exists a
constant $C >
0$ such that for every couple $(t_1, \mathbf{x}_1),  (t_2, \mathbf
{x}_2) \in Q_T$,
%
\begin{eqnarray}
\label{eq:condPhi2} && \bigl|\ell(t_1,\mathbf{x}_1) -
\ell(t_2,\mathbf{x}_2)\bigr| + \bigl| \Phi (\mathbf{x}_1)
- \Phi(\mathbf{x}_2)\bigr |
\nonumber
\\[-8pt]
\\[-8pt]
\nonumber
&&\qquad\le C \bigl( |t_1 - t_2| + \bigl|\mathbf{x}^{t_1}_1
- \mathbf {x}^{t_2}_2\bigr| + |\mathbf{x}_1 -
\mathbf{x}_2| \bigr) \exp \bigl( C\bigl(|\mathbf{x}_1| + |
\mathbf{x}_2| \bigr)\bigr).
\end{eqnarray}

\begin{Theorem} \label{theo:convergence_rate}
Suppose that $L$ and $\Phi$ satisfy conditions \eqref{eq:condPhi1} and
\eqref{eq:condPhi2}, the set $E \subset S_d^+ \times\R^d$ is
compact and
convex, functions $\mu$ and $\sigma$ satisfy \eqref{eq:condUplus}, and
Assumption \ref{assum:G} holds true. Then for every $\eps> 0$, there
is a constant $C_{\eps}$ such that
%
\begin{equation}
\bigl| Y^h_0 - V \bigr| \le C_{\eps} h^{1/8 - \eps}\qquad
\forall h \le h_0.
\end{equation}
If, in addition, $L$ and $\Phi$ are bounded, then there is a constant
$C$ such that
%
\begin{equation}
\bigl| Y^h_0 - V \bigr| \le C h^{1/8}\qquad \forall h \le
h_0.
\end{equation}
\end{Theorem}

\begin{Remark}
In the Markovian context as in Remark \ref{rem:MarkovScheme}, Fahim,
Touzi and Warin \cite{FahimTouziWarin} obtained a convergence rate
$h^{1/4}$ for one side and $h^{1/10}$ for the other side using Krylov's
shaking coefficient method.
Then their global convergence rate is $h^{1/10}$.
We get a rate $h^{1/8}$ in this path-dependent case under some
additional constraints. When there is no control on the volatility
part, the BSDE method in Bouchard and Touzi \cite{Bouchard2004} and
Zhang \cite{Zhang2004} gives a convergence rate of order $h^{1/2}$. Our
current technique cannot achieve this rate in the BSDE context.
\end{Remark}

\begin{Remark}
When the covariance matrix $\sigma\sigma^{\top}$ is diagonal
dominated, Kushner and Dupuis \cite{KushnerDupuis} gave a systematic
way to construct a convergent finite difference scheme. However, the
construction turns to be not easy when the matrix is not diagonal
dominated; see, for example, Bonnans, Ottenwaelter and Zidani~\cite
{Bonnans2004}. Our scheme relaxes this constraint. Moreover, our
scheme implies a natural Monte Carlo implementation, which may be more
efficient in high dimensional cases, see numerical examples in Fahim,
Touzi and Warin \cite{FahimTouziWarin}, Guyon and Henry-Labord\`ere~\cite{GuyonHenryLabordere} and Tan \cite{TanSplitting}.
\end{Remark}

\section{A controlled discrete-time semimartingale interpretation}
\label{sec:proba_interp}

Before giving the proofs of the above convergence theorems, we first
provide a probabilistic interpretation to scheme \eqref{eq:NumSchem},
\eqref{eq:NumTermCond} in the spirit of Kushner and Dupuis \cite
{KushnerDupuis}. Namely, we shall show that the numerical solution is
equivalent to the value function of a controlled discrete-time
semimartingale problem.

For finite difference schemes, the controlled Markov chain
interpretation given by Kushner and Dupuis \cite{KushnerDupuis} is
straightforward, where their construction of the Markov chain is
descriptive. For our scheme, the probabilistic interpretation is less
evident as the state space is uncountable. Our main idea is to use the
inverse function of the distribution functions. This question has not
been evoked in the Markovian context of \cite{FahimTouziWarin} since
they use the monotone convergence of a viscosity solution technique,
where the idea is to show that the terms $Z^h$ and $\Gamma^h$ defined
below \eqref{eq:NumTermCond} are good approximations of the derivatives
of the value function.

\subsection{A technical lemma}
\label{subsec:technique_lemma}

Given a fixed $(t, \mathbf{x}, u) \in Q_T$, let us simplify further the
notation in \eqref{eq:note_a_u},
%
\begin{equation}
\label{eq:note_a_u_simp} \sigma_0:= \sigma_0^{t,\mathbf{x}},\qquad
a_0:= a_0^{t,\mathbf
{x}},\qquad a_u:=
a^{t, \mathbf{x}}_u,\qquad  b_u:= b^{t, \mathbf{x}}_u.
\end{equation}
Denote
%
\begin{eqnarray}
\label{eq:denstiy1} &&f_h(t,\mathbf{x}, u, x)\nonumber \\
&&\qquad:= \frac{1}{(2 \pi h)^{d/2} |\sigma
_0|^{1/2}} \exp
\biggl(-\frac{1}{2} h^{-1} x^{\top} a_0^{-1}
x \biggr)
\\
&&\qquad\quad{}\times \biggl(1 - \frac{1}{2} a_u \cdot a_0^{-1}
+ b_u \cdot a_0^{-1} x + \frac{1}{2}
h^{-1} a_u \cdot a_0^{-1} x
x^{\top} \bigl(a_0^{\top}\bigr)^{-1}
\biggr).\nonumber
\end{eqnarray}
It follows by \eqref{eq:defmG} that for every $(t, \mathbf{x}, u) \in Q_T$
and $x \in\R^d$,
\begin{eqnarray*}
&&b_u \cdot a_0^{-1} x + \frac{1}{2}
h^{-1} a_u \cdot a_0^{-1} x
x^{\top} \bigl(a_0^{\top}\bigr)^{-1} \\
&&\qquad= h
\biggl[ b_u \cdot a_0^{-1} \frac{x}{h} +
\frac{1}{2} a_u \cdot a_0^{-1}
\frac{x}{h} \frac{x^{\top}}{h} \bigl(a_0^{\top}
\bigr)^{-1} \biggr]
\\
&&\qquad\ge h m_G.
\end{eqnarray*}
Then under Assumption \ref{assum:G}, one can verify easily (see also
Remark \ref{rem:positive_density}) that when $h \le h_0$ for $h_0$
given by \eqref{eq:defh0}, $x \mapsto f_h(t,\mathbf{x},u,x)$ is a probability
density function on $\R^d$, that is,
\[
f_h(t,\mathbf{x},u,x) \ge 0\qquad \forall x \in\R^d\quad\mbox
{and}\quad \int_{\R^d} f_h(t,\mathbf{x},u,x) \,dx = 1.
\]

\begin{Lemma} \label{lemma:proba_density}
Let $h \le h_0 $ and $R$ be a random vector with probability density
$x \mapsto f_h(t,\mathbf{x},u,x)$. Then for all functions $g\dvtx \R^d
\to\R$ of
exponential growth, we have
%
\begin{eqnarray}
\label{eq:propR1} && \E\bigl[ g(R) \bigr]= \E \biggl[ g(\sigma_0 W_h) \biggl( 1 + h
b_u \cdot\bigl(\sigma_0^{\top
}
\bigr)^{-1} \frac{W_h}{h}
\nonumber
\\[-8pt]
\\[-8pt]
\nonumber
&&\hspace*{108pt}{}+ \frac{1}{2} h a_u
\cdot\bigl(\sigma_0^{\top}\bigr)^{-1}
\frac{W_h W^{\top}_h - h I}{h^2} \sigma_0^{-1} \biggr) \biggr],
\end{eqnarray}
where $W_h$ is a $d$-dimensional Gaussian random variable with
distribution $N(0, h I_d)$. In particular, it follows that there exists
a constant $C_1$ independent of $(h,t,\mathbf{x}, u) \in(0,h_0]
\times
Q_T$ such that
%
\begin{eqnarray}
\label{eq:propR} \E[R] &=& b_u h, \qquad\operatorname{Var} [R] = (a_u +
a_0) h - b_u b_u^{\top}
h^2 \quad\mbox{and}
\nonumber
\\[-8pt]
\\[-8pt]
\nonumber
\E\bigl[|R|^3\bigr]& <& C_1
h^{3/2},
\end{eqnarray}
where $\operatorname{Var}[R]$ means the covariance matrix of the random vector $Z$.
Moreover, for any $c \in\R^d$,
%
\begin{equation}
\label{eq:propRexp} \E\bigl[ e^{c \cdot R} \bigr] \le e^{C_2 h} (1 +
C_2 h),
\end{equation}
where $C_2$ is independent of $(h,t, \mathbf{x}, u)$ and is defined by
\[
C_2:= \sup_{(t,\mathbf{x}, u) \in Q_T} \biggl( \frac{1}{2}
c^{\top} a_0 c + |b_u \cdot c| +
\frac{1}{2} c^{\top} a_u c \biggr).
\]
\end{Lemma}
\begin{pf} First, it is clear that $\frac{1}{(2\pi h)^{d/2} |\sigma
_0|^{1/2}} \exp (-\frac{1}{2} h^{-1} x^{\top} a_0^{-1} x  )$ is
the density function of $\sigma_0 W_h$. Then by \eqref{eq:denstiy1},
\begin{eqnarray*}
&&\hspace*{-4pt} \E\bigl[g(R)\bigr]\\
 &&\hspace*{-4pt}\qquad = \int_{\R^d} f_h(t,
\mathbf{x}, u,x) g(x) \,dx
\\
&&\hspace*{-4pt}\qquad= \int_{\R^d} \frac{1}{(2\pi h)^{d/2} |\sigma_0|^{1/2}} \exp \biggl(-
\frac{1}{2} h^{-1} x^{\top} a_0^{-1}
x \biggr)
\\
&&\hspace*{-4pt}\hspace*{27pt}\qquad{}\times g(x) \biggl(1 - \frac{1}{2} a_u \cdot a_0^{-1}
+ b_u \cdot a_0^{-1} x + \frac{1}{2}
h^{-1} a_u \cdot a_0^{-1} x
x^{\top} \bigl(a_0^{\top}\bigr)^{-1}
\biggr) \,dx
\\
&&\hspace*{-4pt}\qquad= \E \biggl[ g(\sigma_0 W_h) \biggl( 1 + h
b_u \cdot\bigl(\sigma_0^{\top
}
\bigr)^{-1} \frac{W_h}{h} + \frac{1}{2} h a_u
\cdot\bigl(\sigma_0^{\top}\bigr)^{-1}
\frac{W_h W^{\top}_h - h I}{h^2} \sigma_0^{-1} \biggr) \biggr].
\end{eqnarray*}
Hence \eqref{eq:propR1} holds true.

In particular, let $g(R) = R$ or $g(R) = R R^{\top}$, it follows by
direct computation that the first two equalities of \eqref{eq:propR}
hold true. Further, letting $g(x) = |x|^3$, we get from \eqref
{eq:propR1} that
\begin{eqnarray*}
&&\E\bigl[|R|^3\bigr]\\
&&\qquad = h^{3/2} \E \bigl[ |
\sigma_0 N |^3 \bigl( 1 + \sqrt{h} b_u \cdot
\bigl(\sigma_0^{\top}\bigr)^{-1} N +
\tfrac{1}{2} a_u \cdot\bigl(\sigma _0^{\top}
\bigr)^{-1} \bigl(N N^{\top} - I\bigr) \sigma_0^{-1}
\bigr) \bigr],
\end{eqnarray*}
where $N$ is a Gaussian vector of distribution $N(0, I_d)$. And hence
\eqref{eq:propR} holds true with
\begin{eqnarray*}
&&C_1:= \sup_{(t,\mathbf{x},u) \in Q_T} \E \biggl[ |
\sigma_0 N |^3 \biggl( 1 + \sqrt{h_0}
\bigl|b_u \cdot\bigl(\sigma_0^{\top}
\bigr)^{-1} N \bigr| \\
&&\hspace*{123pt}{}+ \frac{1}{2} \bigl| a_u \cdot\bigl(
\sigma_0^{\top}\bigr)^{-1} \bigl( N N^{\top}
- I \bigr) \sigma _0^{-1} \bigr| \biggr) \biggr],
\end{eqnarray*}
which is clearly bounded and independent of $(h,t,\mathbf{x},u)$.

Finally, to prove inequality \eqref{eq:propRexp}, we denote $N_h:= N
+ \sqrt{h} \sigma_0^{\top} c$ for every $h \le h_0$. Then
\begin{eqnarray*}
&&\E\bigl[ e^{c \cdot R} \bigr]\\
 &&\qquad= \E \biggl[ e^{ c^{\top} \sigma_0 W_h} \biggl( 1 + h
b_u \cdot\bigl(\sigma_0^{\top}
\bigr)^{-1} \frac{W_h}{h} + \frac{1}{2} h a_u
\cdot\bigl(\sigma_0^{\top}\bigr)^{-1}
\frac{W_h W^{\top}_h - h I_d}{h^2} \sigma _0^{-1} \biggr) \biggr]
\\
&&\qquad= \E \biggl[ e^{ c^{\top} \sigma_0 N \sqrt{h} } \biggl( 1 + \sqrt {h} b_u \cdot
\bigl(\sigma_0^{\top}\bigr)^{-1} N +
\frac{1}{2} a_u \cdot\bigl(\sigma _0^{\top
}
\bigr)^{-1} \bigl( N N^{\top} - I_d \bigr)
\sigma_0^{-1} \biggr) \biggr]
\\
&&\qquad= e^{({c^{\top} a_0 c}/{2}) h} \E \biggl[ 1 + \sqrt{h} b_u \cdot \bigl(
\sigma_0^{\top}\bigr)^{-1} N_h +
\frac{1}{2} a_u \cdot\bigl(\sigma_0^{\top
}
\bigr)^{-1} \bigl( N_h N_h^{\top} -
I_d \bigr) \sigma_0^{-1} \biggr]
\\
&&\qquad= e^{({c^{\top} a_0 c}/{2}) h} \biggl( 1 + \biggl(b_u \cdot c +
\frac
{1}{2} c^{\top} a_u c \biggr) h \biggr)
\\
&&\qquad \le e^{C_2 h} (1 + C_2 h),
\end{eqnarray*}
where $C_2:= \sup_{(t,\mathbf{x}, u) \in Q_T}  ( \frac{1}{2}
c^{\top} a_0
c + |b_u \cdot c| + \frac{1}{2} c^{\top} a_u c  )$ is bounded and
independent of $(h, t, \mathbf{x}, u)$.
\end{pf}

\begin{Remark} \label{rem:ab_positive}
Since the random vector $R$ does not degenerate to the Dirac mass, it
follows by \eqref{eq:propR} that under Assumption \ref{assum:G},
\[
\sigma\sigma^{\top}(t,\mathbf{x},u) > \mu\mu^{\top}(t,\mathbf
{x},u) h\qquad \mbox {for every } (t,\mathbf{x},u) \in Q_T, h\le
h_0.
\]
\end{Remark}

With this technical lemma, we can give the proof of Proposition \ref{prop:Y_integ}.

\begin{pf*}{Proof of Proposition \ref{prop:Y_integ}} For the first
assertion, it is enough to prove that $\sup_{\nu\in\Uc} \E[ \exp(C
|X^{\nu}_{\cdot}|) ]$ is bounded by condition \eqref{eq:condPhi1}. Note
that $\mu$ and $\sigma$ are uniformly bounded. When $d = 1$, $X^{\nu}$
is a continuous semimartingale whose finite variation part and
quadratic variation are both bounded by a constant $R_T$ for every $\nu
\in\Uc$. It follows by Dambis--Dubins--Schwarz's time change theorem that
%
\begin{equation}
\label{eq:V_bounded} \sup_{\nu\in\Uc} \E\bigl[ \exp\bigl(C
\bigl|X^{\nu}_{\cdot}\bigr|\bigr) \bigr] \le e^{C
R_T} \E \exp
\Bigl(C \sup_{0 \le t \le R_T} |B_t| \Bigr) < \infty,
\end{equation}
where $B$ is a standard one-dimensional Brownian motion. When $d > 1$,
it is enough to remark that for $X = (X^1, \ldots, X^d)$, $\exp ( C
|X_{\cdot}|  ) \le\exp ( C (|X^1_{\cdot}| + \cdots+
|X^d_{\cdot
}|)  )$; and we then conclude the proof of the first assertion
applying the Cauchy--Schwarz inequality.

We prove the second assertion by backward induction. Given $0 \le k
\le n-1$, $x_0, \ldots, x_k \in\R^d$, we denote by $\widehat{x}$ the
linear interpolation path of $x(t_i):= x_i$; denote also
%
\begin{equation}
\label{eq:defL_k} L_k(x_0, \ldots, x_k,u):=
L(t_k, \widehat{x}_{\cdot},u)\qquad \forall u\in E.
\end{equation}
For the terminal condition, it is clear that $Y^h_n(x_0,\ldots, x_n)$
is of exponential growth in $\max_{0 \le i \le n} |x_n|$ by condition
\eqref{eq:condPhi1}. Now, suppose that
\[
\bigl| Y^h_{k+1}(x_0, \ldots, x_{k+1}) \bigr|
\le C_{k+1} \exp \Bigl( C_{k+1} \max_{0 \le i \le k+1}
|x_i| \Bigr).
\]
Let $R_u$ be a random variable of distribution density $x \mapsto
f_h(t_k,\widehat{x}, u, x)$.
Then it follows by \eqref{eq:NumSchem} and Lemma \ref
{lemma:proba_density} that
%
\begin{eqnarray}
\label{eq:Y_h_PPD} && Y^h_k(x_0, \ldots,
x_k)
\nonumber
\\
&&\qquad= \sup_{u \in E} \biggl\{ h L_k(x_0,
\ldots, x_k,u)\nonumber\\
&&\hspace*{53pt}{} + \E \biggl[ Y^h_{k+1}
(x_0, \ldots, x_k, x_k +
\sigma_0 W_h ),
\nonumber
\\[-8pt]
\\[-8pt]
\nonumber
&&\hspace*{45pt}\quad\qquad{}\times \biggl( 1 + h b^{t_k,\widehat{x}}_u \cdot\bigl(
\sigma_0^{\top
}\bigr)^{-1} \frac{W_h}{h}\\
&&\hspace*{98pt}{} +
\frac{1}{2} h a^{t_k,\widehat{x}}_u \cdot \bigl(
\sigma_0^{\top}\bigr)^{-1} \frac{W_h W_h^{\top} - h I}{h^2}
\sigma_0^{-1} \biggr) \biggr] \biggr\}\nonumber
\\
&&\qquad= \sup_{u \in E} \bigl\{ h L_k(x_0,
\ldots, x_k,u) + \E \bigl[ Y^h_{k+1}
(x_0, \ldots, x_k, x_k + R_u )
\bigr] \bigr\}.
\nonumber
\end{eqnarray}
Therefore by \eqref{eq:condPhi1} and \eqref{eq:propRexp},
\begin{eqnarray*}
&& \bigl| Y^h_k(x_0, \ldots, x_k) \bigr|
\\
&&\qquad\le (C_{k+1} + Ch) \exp \Bigl( (C_{k+1} + Ch) \max
_{0 \le i \le
k}|x_i| \Bigr) \sup_{u \in E}
\E \bigl[ \exp \bigl( C_{k+1} |R_u| \bigr) \bigr]
\\
&&\qquad\le e^{C_2 h} (1 + C_2 h) (C_{k+1} + Ch) \exp
\Bigl( (C_{k+1} + Ch) \max_{0 \le i \le k}|x_i|
\Bigr),
\end{eqnarray*}
where $C$ is the same constant given in \eqref{eq:condPhi1}, and the
constant $C_2$ is from \eqref{eq:propRexp} depending on $C_{k+1}$. We
then conclude the proof.
\end{pf*}

\subsection{The probabilistic interpretation}

In this section, we shall interpret the numerical scheme \eqref
{eq:NumSchem} as the value function of a controlled discrete-time
semimartingale problem. In preparation, let us show how to construct
the random variables with density function $x \mapsto f_h(t,\mathbf{x},u,x)$.
Let $F\dvtx \R\to[0,1]$ be the cumulative distribution function of a
one-dimensional random variable, denote by $F^{-1}\dvtx [0,1] \to\R$ its
generalized inverse function. Then given a random variable $U$ of
uniform distribution $U([0,1])$, it is clear that $F^{-1}(U)$ turns to
be a random variable with distribution $F$. In the multi-dimensional
case, we can convert the problem to the one-dimensional case since $\R
^d$ is isomorphic to $[0,1]$, that is, there is a one-to-one mapping
$\kappa\dvtx  \R^d \to[0,1]$ such that $\kappa$ and $\kappa^{-1}$ are both
Borel measurable; see, for example, Proposition 7.16 and Corollary
7.16.1 of Bertsekas and Shreve \cite{BertsekasShreve}.

Define
\[
F_h(t,\mathbf{x},u,x):= \int_{\kappa(y) \le x}
f_h(t,\mathbf {x},u,y) \kappa(y) \,dy.
\]
It is clear that $x \mapsto F_h(t,\mathbf{x},u,x)$ is the distribution
function of random variable $\kappa(R)$ where $R$ is a random variable
of density function $x \mapsto f_h(t,\mathbf{x},u,x)$. Denote by
$F^{-1}_h(t,\mathbf{x}, u,x)$ the inverse function of $x \mapsto
F_h(t,\mathbf{x},
u,x)$ and
%
\begin{equation}
\label{eq:defHh} H_h(t,\mathbf{x},u,x):= \kappa^{-1}
\bigl( F^{-1}_h(t,\mathbf{x},u,x)\bigr).
\end{equation}
Then given a random variable $U$ of uniform distribution on $[0,1]$,
$F^{-1}_h(t,\mathbf{x},\break u,U)$ has the same distribution of $\kappa(R)$ and
$H_h(t,\mathbf{x},u,U)$ is of distribution density $x \mapsto
f_h(t,\mathbf{x},u,x)$.
In particular, it follows that the expression \eqref{eq:Y_h_PPD} of
numerical solution of scheme \eqref{eq:NumSchem} turns to be
%
\begin{eqnarray}
\label{eq:Y_h_PPD2}\qquad&& Y^h_k(x_0, \ldots,
x_k)
\nonumber
\\[-8pt]
\\[-8pt]
\nonumber
&&\qquad= \sup_{u \in E} \E \bigl[h
L_k(x_0, \ldots, x_k, u)
 + Y^h_{k+1} \bigl(x_0, \ldots,
x_k, x_k + H_h(t_k,
\widehat{x}, u, U) \bigr) \bigr],
\end{eqnarray}
where $\widehat{x}$ is the linear interpolation function of $(x_0,
\ldots, x_k)$ on $[0,t_k]$.

Now, we are ready to introduce a controlled discrete-time
semimartingale system. Suppose that $ U_1, \ldots, U_n$ are i.i.d.
random variables with uniform distribution on $[0,1]$ in the
probability space $(\Om, \Fc, \P)$. Let $\Ac_h$ denote the collection
of all strategies $\phi= (\phi_k)_{0 \le k \le n-1}$, where $\phi_k$
is a universally measurable mapping from $(\R^d)^{k+1}$ to $E$. Given
$\phi\in\Ac_h$, $X^{h,\phi}$ is defined by $X^{h,\phi}_0:= x_0$ and
%
\begin{equation}
\label{eq:defXh}  X^{h,\phi}_{k+1}:= X^{h,\phi}_k
+ H_h \bigl( t_k, \widehat {X}{}^{h,\phi}_{\cdot},
\phi_k\bigl(X^{h,\phi}_0, \ldots,
X^{h,\phi}_k\bigr), U_{k+1} \bigr).
\end{equation}
We then also define an optimization problem by
%
\begin{equation}
\label{eq:Vh0} V^h_0:= \sup_{\phi\in\Ac_h}
\E \Biggl[ \sum_{k=0}^{n-1} h L
\bigl(t_k, \widehat{X}{}^{h,\phi}_{\cdot},
\phi_k\bigr)  + \Phi \bigl(\widehat {X}{}^{h,\phi}_{\cdot}
\bigr) \Biggr].
\end{equation}

The main result of this section is to show that the numerical solution
given by~\eqref{eq:NumSchem} is equivalent to the value function of
optimization problem \eqref{eq:Vh0} on the controlled discrete-time
semimartingales $X^{h,\phi}$.

\begin{Remark}
It is clear that in the discrete-time case, every process is a
semimartingale. When $\mu\equiv0$ and $U$ is of uniform distribution
on $[0,1]$, the random variable $H_h(t,\mathbf{x},u,U)$ is centered,
and hence
$X^{h,\phi}$ turns to be a controlled martingale. This is also the main
reason we choose the terminology ``semimartingale'' in the section title.
\end{Remark}

\begin{Theorem} \label{theo:proba_interp}
Suppose that $L$ and $\Phi$ satisfy \eqref{eq:condPhi1} and Assumption
\ref{assum:G} holds true. Then for $0 < h \le h_0 $ with $h_0$ defined
by \eqref{eq:defh0},
\[
Y^h_0 = V^h_0.
\]
\end{Theorem}

The above theorem is similar to a dynamic programming result. Namely,
it states that optimizing the criteria globally in \eqref{eq:Vh0} is
equivalent to optimizing it step by step in \eqref{eq:Y_h_PPD2}.
With this interpretation, we only need to analyze the ``distance'' of
the controlled semimartingale $X^{h,\phi}$ in \eqref{eq:defXh} and the
controlled diffusion process $X^{\nu}$ in \eqref{eq:SDE} to show this
convergence of $V^h_0$ to $V$ in order to prove Theorems \ref
{theo:convergence} and \ref{theo:convergence_rate}.
Before providing the proof, let us give a technical lemma.
%
\begin{Lemma} \label{Lemma:G_mesurable}
For the function $G$ defined by \eqref{eq:functionG} and every $\eps>
0$, there is a universally measurable mapping $u^{\eps} \dvtx S_d \times
\R^d
\to  E$ such that for all $(\gamma,p) \in S_d \times\R^d$,
\[
G(t, \mathbf{x}, \gamma, p) \le L\bigl(t, \mathbf{x}, u^{\eps}(\gamma,
p)\bigr) + \tfrac
{1}{2} a^{t,\mathbf{x}}_{u^{\eps}(\gamma,p)} \cdot\gamma +
b^{t,\mathbf{x}
}_{u^{\eps}(\gamma,p)} \cdot p + \eps.
\]
\end{Lemma}
\begin{pf}
This follows from the measurable selection
theorem; see, for
example, Theorem 7.50 of Bertsekas and Shreve \cite{BertsekasShreve} or
Section~2 of El Karoui and Tan~\cite{ElKarouiTan}.
\end{pf}

\begin{pf*}{Proof of Theorem \ref{theo:proba_interp}} First,
following \eqref{eq:Y_h_PPD2}, we can rewrite $Y^h_k$ as a measurable
function of $(X^0_0, \ldots, X^0_k)$, and
\begin{eqnarray*}
&& Y^h_k(x_0, \ldots, x_k) \\
&&\qquad=
\sup_{u \in E} \E \bigl[h L_k(x_0, \ldots
, x_k, u) + Y^h_{k+1} \bigl(x_0, \ldots,
x_k, x_k + H_h(t_k,
\widehat{x}, u, U_{k+1}) \bigr) \bigr],
\end{eqnarray*}
where $\widehat{x}$ is the linear interpolation function of $(x_0,
\ldots, x_k)$ on $[0,t_k]$, $U_{k+1}$ is of uniform distribution on
$[0,1]$ and $L_k$ is defined by \eqref{eq:defL_k}.

Next, for every control strategy $\phi\in\Ac_h$ and $X^{h,\phi}$
defined by \eqref{eq:defXh}, we denote $\Fc_k^{h,\phi}:= \sigma(
X^{h,\phi}_0, \ldots, X^{h,\phi}_k)$ and
\[
V^{h,\phi}_k:= \E \Biggl[ \sum
_{i=k}^{n-1} h L\bigl(t_i, \widehat
{X}{}^{h,\phi}_{\cdot}, \phi_i \bigr) + \Phi\bigl(
\widehat{X}{}^{h,\phi
}_{\cdot}\bigr)\Big | \Fc^{h,\phi}_k
\Biggr],
\]
which is clearly a measurable function of $(X^{h,\phi}_0, \ldots,
X^{h,\phi}_k)$ and satisfies
\begin{eqnarray*}
V^{h,\phi}_k (x_0,\ldots,x_k) &=& h
L_k\bigl( x_0, \ldots, x_k,
\phi_k(x_0, \ldots, x_k)\bigr)
\\
&&{} + \E \bigl[ V^{h,\phi}_{k+1} \bigl(x_0, \ldots,
x_k, x_k + H_h\bigl(t_k,
\widehat{x}, \phi_k(x_0,\ldots, x_k),
U_{k+1}\bigr) \bigr) \bigr].
\end{eqnarray*}
Then by comparing $V^{h,\phi}$ with $Y^h$ and the arbitrariness of
$\phi\in\Ac_h$, it follows that
\[
V^h_0 \le Y^h_0.
\]

For the reverse inequality, it is enough to find, for any $\eps> 0$,
a strategy $\phi^{\eps} \in\Ac_h$ with $X^{h,\eps}$ as defined in
\eqref{eq:defXh} using $\phi^{\eps}$ such that
%
\begin{equation}
\label{eq:Xh_eps} Y^h_0 \le \E \Biggl[ \sum
_{k=0}^{n-1} h L\bigl(t_k,
\widehat{X}{}^{h,\eps
}_{\cdot}, \phi^{\eps}_k
\bigr) + \Phi\bigl( \widehat{X}{}^{h,\eps}_{\cdot
}\bigr) \Biggr] + n
\eps.
\end{equation}
Let us write $\Gamma_k^h$ and $Z^h_k$ defined below \eqref
{eq:NumSchem} as a measurable function of $(X_0^0, \ldots, X_k^0)$, and
$u^{\eps}$ be given by Lemma \ref{Lemma:G_mesurable}, denote
\[
\phi^{\eps}_k(x_0, \ldots, x_k)
:= u^{\eps}\bigl(\Gamma^h_k(x_0,
\ldots, x_k), Z^h_k(x_0, \ldots,
x_kz)\bigr).
\]
Then by the tower property, the semimartingale $X^{h,\eps}$ defined by
\eqref{eq:defXh} with $\phi^{\eps}$ satisfies \eqref{eq:Xh_eps}.
\end{pf*}

\section{Proofs of the convergence theorems}
\label{sec:conv_proof}

With the probabilistic interpretation of the numerical solution $Y^h$
in Theorem \ref{theo:proba_interp}, we are ready to give the proofs of
Theorems \ref{theo:convergence} and \ref{theo:convergence_rate}.
Intuitively, we shall analyze the ``convergence'' of the controlled
semimartingale $X^{h,\phi}$ in \eqref{eq:defXh} to the controlled
diffusion process $X^{\nu}$ in~\eqref{eq:SDE}.

\subsection{Proof of Theorem \texorpdfstring{\protect\ref{theo:convergence}}{2.7}}

The main tool we use to prove Theorem \ref{theo:convergence} is the
weak convergence technique due to Kushner and Dupuis \cite
{KushnerDupuis}. We adapt their idea in our context. We shall also
introduce an enlarged canonical space for control problems following El
Karoui, H{\.u}{\.u} Nguyen and Jeanblanc \cite{ElKaroui1987}, in order to
explore the convergence conditions. Then we study the weak convergence
of probability measures on the enlarged canonical space.

\subsubsection{An enlarged canonical space}

In Dolinsky, Nutz and Soner \cite{DolinskyNutzSoner}, the authors
studied a similar but simpler problem in the context of
$G$-expectation, where they use the canonical space $\Om^d:= C([0,T],
\R^d)$. We refer to Stroock and Varadhan \cite{Stroock1979} for a
presentation of basic properties of canonical space $\Om^d$. However,
we shall use an enlarged canonical space introduced by El Karoui, H{\.u}{\.u}
Nguyen and Jeanblanc \cite{ElKaroui1987}, which is more convenient to
study the control problem for the purpose of numerical analysis.

\paragraph*{An enlarged canonical space}
Let $\Mbf([0,T] \times E)$ denote the space of all finite positive
measures $m$ on $[0,T] \times E$ such that $m([0,T] \times E) = T$,
which is a
Polish space equipped with the weak convergence topology. Denote by
$\Mbf$ the collection of finite positive measures $m \in\Mbf
([0,T]\times
E)$ such that the projection of $m$ on $[0,T]$ is the Lebesgue measure,
so that they admit the disintegration $m(dt,du) = m(t,du) \,dt$, where
$m(t,du)$ is a probability measure on $E$ for every $t \in[0,T]$, that is,
\[
\Mbf:= \biggl\{ m \in\Mbf\bigl([0,T]\times E\bigr) \dvtx m(dt, du) = m(t,du) \,dt
\mbox{ s.t. } \int_E m(t,du) = 1 \biggr\}.
\]
In particular, $(m(t,du))_{0 \le t \le T}$ is a measure-valued
process. The measures in space $\Mbf$ are examples of Young measures
and have been largely used in deterministic control problems. We also
refer to Young \cite{Young} and Valadier \cite{Valadier} for a
presentation of Young measure as well as its applications.

Clearly, $\Mbf$ is closed under weak convergence topology and hence is
also a Polish space. We define also the $\sigma$-fields on $\Mbf$ by
$\Mc_t:= \sigma\{ m_s(\varphi),   s\le t,   \varphi\in C_b([0,T]
\times
E) \}$, where $m_s(\varphi):= \int_0^s \varphi(r,u) m(dr, du)$. Then
$(\Mc_t)_{0 \le t \le T}$ turns to be a filtration. In particular,
$\Mc
_T$ is the Borel $\sigma$-field of $\Mbf$. As defined above, $\Om^d:=
C([0,T], \R^d)$ is the space of all continuous paths between $0$ and
$T$ equipped with canonical filtration $\F^d = (\Fc^d_t)_{0 \le t \le
T}$. We then define an enlarged canonical space by $\Omb^d:= \Om^d
\times
\Mbf$, as well as the canonical process $X$ by $X_t(\overline{\omega
}):= \omega^d_t$,
$\forall\overline{\omega}= (\omega^d,m) \in\Omb^d = \Om^d \times
\Mbf$, and the
canonical filtration $\Fbb^d = (\Fcb^d_t)_{0 \le t \le T}$ with $\Fcb
^d_t:= \Fc^d_t \otimes\Mc_t$. Denote also by $\Mbf(\Omb^d)$ the
collection of all probability measures on $\Omb^d$.

\paragraph*{Four classes of probability measures}
A controlled diffusion process as well as the control process may
induce a probability measure on $\Omb^d$. Further, the optimization
criterion in \eqref{eq:V} can be then given as a random variable
defined on $\Omb^d$. Then the optimization problem \eqref{eq:V} can be
studied on $\Omb^d$, as the quasi-sure approach in Soner, Touzi and
Zhang \cite{STZquasisure}. In the following, we introduce four
subclasses of $\Mbf(\Omb^d)$.

Let $\delta> 0$,\vspace*{-1.5pt} we consider a particular strategy $\nu^{\delta}$ of
the form $\nu^{\delta}_s = w_k(X^{\nu^{\delta}}_{r_i^k},\break   i\le I_k)$
for every $s \in(k \delta, (k+1) \delta]$, where $I_k \in\N$, $0
\le
r_0^k < \cdots< r_{I_k}^k \le k \delta$, $X^{\nu^{\delta}}$ is the
controlled process given by \eqref{eq:SDE} with strategy $\nu^{\delta
}$, and $w_k \dvtx \R^{d I_k} \to E$ is a continuous function. Clearly,
$\nu
^{\delta}$ is an adapted piecewise constant strategy. Denote by $\Uc_0$
the collection of all strategies of this form for all $\delta> 0$. It
is clear that $\Uc_0 \subset\Uc$.

Given $\nu\in\Uc$ in the probability space $(\Om, \Fc, \P)$, denote
$m^{\nu}(dt, du):=\break   \delta_{\nu_t}(du) \,dt \in\Mbf$. Then $(X^{\nu},
m^{\nu})$ can induce a probability measure $\Pb^{\nu}$ on $\Omb^d$ by
%
\begin{equation}
\label{eq:defPb} \E^{\Pb^{\nu}} \Upsilon \bigl(\omega^d, m \bigr)
:= \E^{\P} \Upsilon\bigl( X^{\nu}_{\cdot},
m^{\nu} \bigr)
\end{equation}
for every bounded measurable function $\Upsilon$ defined on $\Omb^d$.
In particular, for any bounded function $f \dvtx \R^{dI + IJ} \to\R$ with
arbitrary $I,J \in\N$, $s_i \in[0,T]$, $\psi_j \dvtx [0,T] \times E
\to\R
$ bounded,
\begin{eqnarray*}
&&\E^{\Pb^{\nu}} f \bigl(X_{s_i}, m_{s_i}(
\psi_j), i\le I, j\le J \bigr) \\
&&\qquad= \E^{\P} f
\biggl(X^{\nu}_{s_i}, \int_0^{s_i}
\psi_j(\nu_r) \,dr, i\le I, j\le J \biggr).
\end{eqnarray*}
Then the first and the second subsets of $\Mbf(\Omb^d)$ are given by
\[
\Pcb_{S_0}:= \bigl\{ \Pb^{\nu} \dvtx \nu\in\Uc_0
\bigr\} \quad\mbox {and}\quad \Pcb_S:= \bigl\{ \Pb^{\nu} \dvtx
\nu\in\Uc \bigr\}.
\]

Now, let $ 0 < h \le h_0$ and $\phi\in\Ac_h$, denote $m^{h,\phi}
(dt, du):= \sum_{k=0}^{n-1} \delta_{\phi_k}(du)\times
1_{(t_k,t_{k+1}]}(dt)$, and $X^{h, \phi}$ be the discrete-time
semimartingale defined by \eqref{eq:defXh}. It follows that $(\widehat
{X}{}^{h,\phi}, m^{h,\phi})$ induces a probability measure $\Pb
^{h,\phi}$
on $\Omb^d$ as in \eqref{eq:defPb}. Then the third subset of $\Mbf
(\Omb
^d)$ we introduce is
\[
\Pcb_h:=\bigl\{ \Pb^{h,\phi} \dvtx \phi\in
\Ac_h \bigr\}.
\]

Finally, for the fourth subset, we introduce a martingale problem on
$\Omb^d$. Let $\Lc^{t,\mathbf{x},u}$ be a functional operator
defined by
%
\begin{equation}
\label{eq:defLc} \Lc^{t,\mathbf{x},u} \varphi:= \mu(t,\mathbf{x},u) \cdot D\varphi
+ \tfrac
{1}{2} \sigma\sigma^{\top}(t,\mathbf{x},u) \cdot
D^2 \varphi.
\end{equation}
Then for every $\varphi\in C_b^{\infty}(\R^d)$, a process $M(\varphi
)$ is defined on $\Omb^d$ by
%
\begin{equation}
\label{eq:defM_phi} M_t(\varphi):= \varphi(X_t) -
\varphi(X_0) - \int_0^t \int
_E \Lc ^{t,X_{\cdot},u} \varphi(X_s) m( s, du)
\,d s.
\end{equation}
Denote by $\Pcb_R$ the collection of all probability measures on $\Omb
^d$ under which $X_0 = x_0$ a.s. and $M_t(\varphi)$ is a $\Fbb
^d$-martingale for every $\varphi\in C_b^{\infty}(\R^d)$. In \cite
{ElKaroui1987}, a probability measure in $\Pcb_R$ is called a relaxed
control rule.

\begin{Remark}\label{rem:RelaxMartProb}
We denote, by abuse of notation, the random processes in $\Omb^d$
\begin{eqnarray*}
\mu\bigl(t,\omega^d, m\bigr)&:=& \int_E \mu
\bigl(t,\omega^d,u\bigr) m(t,du),\\
 a\bigl(t,\omega^d, m
\bigr)&:=& \int_E \sigma\sigma^{\top}\bigl(t,
\omega^d,u\bigr) m(t,du),
\end{eqnarray*}
which are clearly adapted to the filtration $\Fbb^d$.
It follows that $M_t(\varphi)$ defined by \eqref{eq:defM_phi} is
equivalent to
\begin{eqnarray*}
M_t(\varphi) &:=& \varphi(X_t) - \varphi(X_0)
\\
&&{}- \int_0^t \biggl( \mu \bigl(s,\omega
^d,m\bigr) \cdot D \varphi(X_s) + \frac{1}{2} a
\bigl(s,\omega^d, m\bigr) \cdot D^2 \varphi
(X_s) \biggr)\,ds.
\end{eqnarray*}
Therefore, under any probability $\Pb\in\Pcb_R$, since $a(s,\omega^d,
m)$ is nondegenerate, there is a Brownian motion $\tilde W$ on $\Omb
^d$ such that the canonical process can be represented as
\[
X_t = x_0 + \int_0^t
\mu(s, X_{\cdot}, m) \,d s + \int_0^t
a^{1/2}(s, X_{\cdot}, m) \,d \tilde W_s.
\]
\end{Remark}

Moreover, it follows by It\^o's formula as well as the definition of
$\Pb^{\nu}$ in \eqref{eq:defPb} that $\Pb^{\nu} \in\Pcb_R$ for every
$\nu\in\Uc$. In resume, we have $\Pcb_{S_0} \subset\Pcb_S \subset
\Pcb_R$.

\paragraph*{A completion space of $C_b(\Omb^d)$}

Now, let us introduce two random variables on $\Omb^d$ by
%
\begin{equation}
\label{eq:defPsi} \Psi(\overline{\omega}) = \Psi\bigl(\omega^d, m\bigr)
:= \int_0^T L \bigl(t,\omega^d, u
\bigr) m(dt,du) + \Phi\bigl(\omega^d\bigr)
\end{equation}
and
\[
\Psi_h(\overline{\omega}) = \Psi_h\bigl(
\omega^d, m\bigr):= \int_0^T
L_h \bigl(t,\omega^d, u \bigr) m(dt,du) + \Phi\bigl(
\omega^d\bigr),
\]
where $L_h(t,\mathbf{x},u):= L(t_k, \mathbf{x}, u)$ for every $t_k
\le t \le
t_{k+1}$ given the discretization parameter $h:= \frac{T}{n}$. It
follows by the uniform continuity of $L$ that
%
\begin{equation}
\label{eq:err_Psi_Psih} \sup_{\overline{\omega}\in\Omb} \bigl|\Psi(\overline{\omega}) -
\Psi_h(\overline{\omega})\bigr| \le \rho_0(h),
\end{equation}
where $\rho_0$ is the continuity module of $L$ in $t$ given before
\eqref{eq:condPhi1}.
Moreover, optimization problems \eqref{eq:V} and \eqref{eq:Vh0} are
equivalent to
%
\begin{equation}
\label{eq:equiv_V_Vh_relax} V = \sup_{\Pb\in\Pcb_S} \E^{\Pb} [ \Psi]
\quad\mbox{and}\quad V_0^h = \sup_{\Pb\in\Pcb_h}
\E^{\Pb} [ \Psi_h ].
\end{equation}

Finally, we introduce a space $L^1_*$ of random variables on $\Omb^d$. Let
\[
\Pcb^*:= \biggl( \bigcup_{0 < h \le h_0}
\Pcb_h \biggr) \cup \Pcb_R,
\]
and defined a norm $|\cdot|_*$ for random variables on $\Omb^d$ by
\[
| \xi|_*:= \sup_{\Pb\in\Pcb^*} \E^{\Pb} | \xi|.
\]
Denote by $L^1_*$ the completion space of $C_b(\Omb^d)$ under the norm
$| \cdot|_*$.

\subsubsection{Convergence in the enlarged space}

We first give a convergence result for random variables in $L^1_*$.
Then we show that $\Psi$ defined by \eqref{eq:defPsi} belongs to
$L^1_*$. In the end, we provide two other convergence lemmas.

\begin{Lemma} \label{lemma:L1starConv}
Suppose that $\xi\in L^1_*$, $(\Pb_n)_{n \ge0}$ is a sequence of
probability measures in $\Pcb^*$ such that $\Pb_n$ converges weakly to
$\Pb\in\Pcb^*$. Then
%
\begin{equation}
\label{eq:weak_conv_L1} \E^{\Pb_n} [\xi] \to \E^{\Pb} [\xi].
\end{equation}
\end{Lemma}
\begin{pf}
For every $\eps> 0$, there is $\xi_{\eps} \in
C_b(\Omb^d)$ such
that $\sup_{\Pb\in\Pcb^*} \E^{\Pb} [| \xi- \xi_{\eps} |]  \le
\eps$.
It follows that
\begin{eqnarray*}
&& \limsup_{n \to\infty} \bigl| \E^{\Pb_n} [\xi] -
\E^{\Pb} [\xi] \bigr|
\\
&&\qquad\le \limsup_{n \to\infty} \bigl[ \E^{\Pb_n} \bigl[ | \xi - \xi
_{\eps} | \bigr] + \bigl| \E^{\Pb_n} [\xi_{\eps}] - \E
^{\Pb} [\xi _{\eps}]\bigr | + \E^{\Pb}\bigl [ |
\xi_{\eps} - \xi | \bigr] \bigr]
\\
&&\qquad\le 2 \eps.
\end{eqnarray*}
Therefore, \eqref{eq:weak_conv_L1} holds true by the arbitrariness of
$\eps$.
\end{pf}

The next result shows that the random variable $\Psi$ defined by
\eqref
{eq:defPsi} belongs to~$L^1_*$, when $L$ and $\Phi$ satisfy \eqref
{eq:condPhi1}.

\begin{Lemma} \label{lemma:L1star}
Suppose that Assumption \ref{assum:G} holds true, and $\Phi$ and $L$
satisfy~\eqref{eq:condPhi1}. Then the random variable $\Psi$ defined by
\eqref{eq:defPsi} lies in $L^1_*$.
\end{Lemma}
\begin{pf}
We first claim that for every $C > 0$
%
\begin{equation}
\label{eq:claim_exp_bound} \sup_{\Pb\in\Pcb^*} \E^{\Pb} \bigl[
e^{C |X_T|} \bigr] < \infty,
\end{equation}
which implies that $\E[e^{C |X^{h,\phi}_n|}] < \infty$ is uniformly
bounded in $h$ and in $\phi\in\Ac_h$. Let $C_0 > 0$ such that $|\mu
(t,\mathbf{x},u)| \le C_0$, $\forall(t,\mathbf{x},u) \in Q_T$. Then
$(|X^{h,\phi}_k
- C_0 t_k|)_{0 \le k \le n}$ is a submartingale, and hence for every
$C_1 > 0$,\break  $( e^{ C_1 |X^{h,\phi}_k - C_0 t_k|})_{0 \le k \le n}$ is
also a submartingale. Therefore, by Doob's inequality,
\begin{eqnarray*}
\E \Bigl[ \sup_{0 \le k \le n} e^{C_1 |X^{h,\phi}_k|} \Bigr] &\le&
e^{d
C_0 T} \E \Bigl[ \sup_{0 \le k \le n} e^{C_1 |X^{h,\phi}_k - C_0 t_k |} \Bigr]
\\
&\le& 2 e^{2 d C_0 T} \sqrt{ \E \bigl[ e^{2C_1 |X^{h,\phi}_n|} \bigr]},
\end{eqnarray*}
where the last term is bounded uniformly in $h$ and $\phi$ by the
claim \eqref{eq:claim_exp_bound}. With the same arguments for the
continuous-time case in spirit of Remark \ref{rem:RelaxMartProb}, it
follows by~\eqref{eq:condPhi1} and \eqref{eq:defPsi} that
%
\begin{equation}
\label{eq:Psi2_bounded} \sup_{\Pb\in\Pcb^*} \E^{\Pb} \bigl[ |
\Psi|^2 \bigr]\le \infty.
\end{equation}
Similarly, we also have $\sup_{\Pb\in\Pcb^*}  \E^{\Pb}    [
|\Psi
'|^2  ] \le\infty$ for
\[
\Psi'\bigl(\omega^d, m\bigr):= \int
_0^T \bigl| L \bigl(t,\omega^d, u\bigr) \bigr|
m(dt,du) + \bigl|\Phi\bigl(\omega^d\bigr) \bigr|.
\]

Let $\Phi_N:= (-N) \vee(\Phi\land N)$, $ L_N:= (-N) \vee(L \land
N)$ and
\[
\Psi_N\bigl(\omega^d, m\bigr):= \int
_0^T \int_E
L_N \bigl(t,\omega^d, u\bigr) m(dt,du) +
\Phi_N\bigl(\omega^d\bigr).
\]
Then $\Psi_N$ is bounded continuous in $ \overline{\omega}= (
\omega^d, m)$, that is, $
\Psi_N \in C_b(\Omb^d)$. It follows by the Cauchy--Schwarz inequality that
\begin{eqnarray*}
\sup_{\Pb\in\Pcb^*} \E^{\Pb} |\Psi- \Psi_N| &
\le& \sqrt{ \sup_{\Pb\in\Pcb^*} \E^{\Pb} |\Psi-
\Psi_N|^2 } \sqrt{ \sup_{\Pb
\in
\Pcb^*} \Pb
\bigl(\bigl|\Psi'\bigr| > N\bigr) }
\\
&\le& \sqrt{ \sup_{\Pb\in\Pcb^*} \E^{\Pb} |\Psi-
\Psi_N|^2 } \sqrt { \sup_{\Pb\in\Pcb^*}
\E^{\Pb} \bigl[ \bigl|\Psi'\bigr| \bigr] } \frac{1}{\sqrt{N}} \to 0,
\end{eqnarray*}
where the last inequality is from $\Pb(|\Psi'| > N)  \le  \frac{1}{N}
\E^{\Pb} [|\Psi'|] $.
And hence $\Psi\in L^1_*$. Therefore, it is enough to justify claim
\eqref{eq:claim_exp_bound} to complete the proof.

By Lemma \ref{lemma:proba_density}, for every random variable $R$ of
density function $f_h(t,\mathbf{x},u,x)$ and every $c \in\R^d$, we have
\[
\E\bigl[ e^{c \cdot R} \bigr] \le e^{C_2 h} (1 + C_2
h),
\]
where $C_2:= \sup_{(t,\mathbf{x}, u) \in Q_T}  ( \frac{1}{2}
|c^{\top}
a^{t,\mathbf{x}}_0 c| + |b^{t,\mathbf{x}}_u \cdot c| + \frac{1}{2}
|c^{\top} a^{t,\mathbf{x}
}_u c|  )$. It follows by taking conditional expectation on $e^{c
\cdot X^h_n}$ that
\[
\E\bigl[ e^{c \cdot X^h_n} \bigr] \le C_0(c):= \sup
_{h \le h_0} e^{C_2 T} (1+ C_2h)^{T/h} <
\infty.
\]
Let $c$ be the vectors of the form $(0, \ldots, 0, \pm C,0, \ldots,
0)^{\top}$, and we can easily conclude that $\E e^{C|X^h_n|}$ is
uniformly bounded for all $h\le h_0$ and $X^h = X^{h,\phi}$ with $\phi
\in\Ac_h$. Furthermore, in spirit of Remark \ref{rem:RelaxMartProb}
and by the same arguments as \eqref{eq:V_bounded} in the proof of
Proposition \ref{prop:Y_integ}, $\sup_{\Pb\in\Pcb_R} \E^{\Pb} [ e^{C
| X_T| } ]$ is bounded. And therefore, we proved the claim \eqref
{eq:claim_exp_bound}.
\end{pf}

Finally, we finish this section by providing two convergence lemmas,
but leave their proofs in \hyperref[app]{Appendix}.

\begin{Lemma} \label{lemma:weakCvg}
\textup{(i)} Let $(\Pb_h)_{0 < h \le h_0}$ be a sequence of probability
measures such that $\Pb_h \in\Pcb_h$. Then $(\Pb_h)_{0 < h \le h_0}$
is precompact, and any cluster point belongs to $\Pcb_R$.

\textup{(ii)} Let $\Pb\in\Pcb_{S_0}$. Then we can construct a
sequence of
probability measures $(\Pb_h)_{0 < h \le h_0}$ such that $\Pb_h \in
\Pcb
_h$ and $\Pb_h \to\Pb$ as $h \to0$.
\end{Lemma}

\begin{Lemma} \label{lemma:equiv_WS}
Suppose that Assumptions \ref{assum:G} holds true, and $\Phi$ and $L$
satisfy~\eqref{eq:condPhi1}. Then
\[
\sup_{\Pb\in\Pcb_{S_0}} \E^{\Pb} [\Psi] = \sup
_{\Pb\in\Pcb
_R} \E ^{\Pb} [\Psi].
\]
\end{Lemma}

\subsubsection{Proof of the general convergence (Theorem \texorpdfstring{\protect\ref{theo:convergence}}{2.7})}

Finally, we are ready to give the proof of Theorem \ref{theo:convergence}.

\begin{pf*}{Proof of Theorem \ref{theo:convergence}} Since $\Psi
\in L^*_1$ by Lemma \ref{lemma:L1star}, then in spirit of Lemma \ref
{lemma:L1starConv}, we get from (i) of Lemma \ref
{lemma:weakCvg} that
\[
\limsup_{h \to0} \sup_{\Pb_h \in\Pcb_h} \E^{\Pb_h}
[\Psi] \le  \sup_{\Pb\in\Pcb_R} \E^{\Pb} [\Psi].
\]
Moreover, it follows by (ii) of Lemma \ref{lemma:weakCvg} that
\[
\liminf_{h \to0} \sup_{\Pb_h \in\Pcb_h} \E^{\Pb_h}
[\Psi] \ge  \sup_{\Pb\in\Pcb_{S_0}} \E^{\Pb} [\Psi].
\]
We hence conclude the proof of the theorem by Lemma \ref
{lemma:equiv_WS} and \eqref{eq:err_Psi_Psih}, \eqref
{eq:equiv_V_Vh_relax}.
\end{pf*}

\subsection{Proofs of Theorem \texorpdfstring{\protect\ref{theo:convergence_rate}}{2.8}}

The proof of Theorem \ref{theo:convergence_rate} is similar to
Dolinsky~\cite{Dolinsky}, where the author uses the invariance
principle technique of Sakhanenko~\cite{Sakhanenko} to approximate the
discrete-time martingales. In our context, we shall approximate
discrete-time semimartingales.

\subsubsection{From continuous to discrete-time semimartingale}
\label{subsubsec:disc_semimart}

The next result is similar to Lemmas 4.2 and 4.3 of Dolinsky \cite
{Dolinsky}, which states that a continuous martingale can be
approximated by its discrete-time version.\vadjust{\goodbreak}

In \eqref{eq:condUplus}, we assume that $E$ is a convex compact subset
in $S^+_d \times\R^d$ and $\mu(t,\mathbf{x},u) = b$, $\sigma
(t,\mathbf
{x},u) = a^{1/2}$
for every $u = (a,b) \in E$. Given a strategy $\nu= (a_t, b_t)_{0 \le
t \le T} \in\Uc$ as well as a discrete time grid $\pi= (t_k)_{0 \le k
\le n}$ ($t_k:= kh$, $h:= T/n \le h_0$), let us define the following
discrete-time processes:
%
\begin{eqnarray}
B_0^{\pi,\nu}:= 0,\qquad B_{k+1}^{\pi,\nu}:=
B_k^{\pi,\nu} + \Delta B_{k+1}^{\pi,\nu}
\nonumber\\
\eqntext{\mbox{with }\displaystyle\Delta B_{k+1}^{\pi,\nu}:= \E _k^{\pi}
\biggl[ \int_{t_k}^{t_{k+1}} b_s\,ds
\biggr],}
\end{eqnarray}
%
\begin{equation}
\label{eq:defDeltaM} \Delta M^{\pi,\nu}_{k+1}:= \int
_{t_k}^{t_{k+1}} \bigl(a^b_s
\bigr)^{1/2} \,dW_s
\end{equation}
with $a^b_s:=  a_s  -  \Delta B_{k+1}^{\pi,\nu}  (\Delta
B_{k+1}^{\pi,\nu}  )^{\top}/h$,
\[
M^{\pi,\nu}_{k+1}:= M^{\pi,\nu}_k + \Delta
M^{\pi,\nu}_{k+1} \quad\mbox{and}\quad X_k^{\pi,\nu}:=
x_0 + B^{\pi,\nu}_k + M^{\pi
,\nu}_k,
\]
where $\E^\pi_k [\cdot]:= \E[ \cdot| \Fc^\pi_k]$ for $\Fc^\pi
_k:=
\sigma(X_0^{\pi,\nu}, \ldots, X_k^{\pi,\nu})$.
We notice that by Remark~\ref{rem:ab_positive}, the matrix $a_s^b$
defined by \eqref{eq:defDeltaM} is strictly positive for every $s$
under Assumption \ref{assum:G} and hence $\Delta M^{\pi,\nu}$ is
well defined.
We denote also
\[
\Delta A^{\pi,\nu}_{k+1}:= \E_k^\pi
\bigl[ \Delta M^{\pi,\nu}_{k+1} \bigl(M^{\pi,\nu}_{k+1}
\bigr)^{\top} \bigr] = \E_k^\pi \biggl[ \int
_{t_k}^{t_{k+1}} a_s \,ds \biggr] - \Delta
B_{k+1}^{\pi,\nu} \bigl(\Delta B_{k+1}^{\pi,\nu}
\bigr)^{\top}.
\]
Similarly, for every $(x, \nu) = (x_0, \ldots, x_n, \nu_1, \ldots,
\nu
_n) \in\R^{d (n+1)} \times(S_d^+ \times\R^d)^n$, a~discrete-time
version of
function $\Psi$ in \eqref{eq:defPsi} can be given by
%
\begin{equation}
\label{eq:defPsi_h} \Psi_\pi(x,\nu):= \sum
_{k=0}^{n-1} h L ( t_k, \widehat
{x}_{\cdot
}, \nu_{k+1} ) + \Phi(\widehat{x}_{\cdot}),
\end{equation}
where $\hat{x}$ is the linear interpolation function of $x_0, \ldots, x_n$.

Now we introduce a discrete-time version of optimization problem \eqref{eq:V},
%
\begin{equation}\qquad
\label{eq:Vd} V^\pi:= \sup_{\nu\in\Uc} \E \bigl[
\Psi_{\pi}\bigl( X^{\pi
,\nu}, \nu^{\pi,\nu}\bigr) \bigr]\qquad
\mbox{with } \nu^{\pi,\nu}_k:= \biggl( \frac
{1}{h} \Delta
A^{\pi,\nu}_k, \frac{1}{h} \Delta B^{\pi,\nu}_k
\biggr).
\end{equation}

\begin{Remark}
The definition of $M^{\pi,\nu}$ in \eqref{eq:defDeltaM} uses a
perturbation version of~$a$. The main purpose is to adapt the biased
term appearing in the variance term of~\eqref{eq:propR}. In particular,
it follows that $\nu^{\pi,\nu}_k \in E_h$ almost surely for
%
\begin{equation}
\label{eq:defEh} E_h:= \bigl\{\bigl(a - hb b^{\top},b\bigr)
\dvtx (a,b) \in E \bigr\}.
\end{equation}
\end{Remark}

\begin{Lemma} \label{lemma:error_discretSemiMart}
\textup{(i)} There is a constant $C$ independent of $h = 1/n$ such that
%
\begin{equation}
\label{eq:error_discretSemiMart} \sup_{\nu\in\Uc} \E \bigl|X^{\nu}_{\cdot}
- \widehat{X}{}^{\pi
,\nu
}_{\cdot} \bigr|^2 \le C
h^{1/2}.
\end{equation}
\textup{(ii)}
It follows that under conditions \eqref{eq:condPhi1}
and \eqref
{eq:condPhi2}, there is a constant $C$ such that
%
\begin{equation}
\label{eq:error_discretV} \bigl| V - V^\pi \bigr| \le C h^{1/4}.
\end{equation}
\end{Lemma}
\begin{pf}
\textup{(i)} Given a control $\nu= (a_t,
b_t)_{0 \le t \le T}
\in\Uc
$, denote
\[
B^{\nu}_t:= \int_0^t
b_s \,d s,\qquad \tilde M_t^{\nu}:= \int
_0^t a^{1/2}_s
\,dW_s \quad\mbox{and}\quad M_t^{\nu}:= \int
_0^t \bigl(a^b_s
\bigr)^{1/2} \,dW_s.
\]
Then it is clear that $X^{\nu}_t = X_0 + B^{\nu}_t + \tilde M_t^{\nu
}$. Since $\nu_s = (a_s, b_s)$ is uniformly bounded, there is a
constant $C$ independent of $n$ such that
%
\begin{equation}
\sup_{\nu\in\Uc} \E \biggl| \int_0^T
\nu_s \,ds - \sum_{k=1}^n h
\nu ^{\pi,\nu}_k \biggr|^2 \le C \frac{1}{n}
\end{equation}
and
\[
\sup_{\nu\in\Uc} \E \bigl|B^{\nu}_{\cdot} -
\widehat{B}{}^{\pi
,\nu
}_{\cdot}\bigr |^2 \le C
\frac{1}{n}.
\]
Moreover,
it follows by Lemmas 4.2 of Dolinsky \cite{Dolinsky} that
\[
\sup_{\nu\in\Uc} \E \bigl| \tilde M^{\nu}_{\cdot} -
\widehat {M}{}^{\pi
,\nu}_{\cdot}\bigr |^2 \le 2 \sup
_{\nu\in\Uc} \E \bigl[\bigl | \tilde M^{\nu}_{\cdot} -
M^{\nu}_{\cdot}\bigr |^2 + \E \bigl| M^{\nu
}_{\cdot}
- \widehat{M}{}^{\pi,\nu}_{\cdot} \bigr|^2 \bigr]
 \le C
\frac
{1}{\sqrt{n}}.
\]
Therefore, by the fact that $|X^{\nu}_{\cdot} - \widehat{X}{}^{\pi,\nu
}_{\cdot} |^2  \le  2  ( |B^{\nu}_{\cdot} - \widehat{B}{}^{\pi
,\nu
}_{\cdot} |^2 + | \tilde M^{\nu}_{\cdot} - \widehat{M}{}^{\pi,\nu
}_{\cdot
} |^2  )$, we prove \eqref{eq:error_discretSemiMart}.

(ii) For the second assertion, we remark that by \eqref
{eq:condPhi2},
for every $\nu\in\Uc$,
\begin{eqnarray*}
&&\Biggl | \int_0^T L\bigl(s,X^{\nu}_{\cdot},
\nu_s\bigr) \,ds - \sum_{k=0}^{n-1}
h L \bigl( t_k, \widehat{X}{}^{\pi,\nu}_{\cdot},
\nu^{\pi,\nu}_{k+1} \bigr) + \Phi\bigl(X^{\nu}_{\cdot}
\bigr) - \Phi\bigl(\widehat{X}{}^{\pi,\nu}_{\cdot}\bigr)\Biggr |
\\
&&\qquad\le C \exp \bigl(C\bigl(\bigl|X^{\nu}_{\cdot}\bigr| + \bigl|
\widehat{X}{}^{\pi,\nu
}_{\cdot
}\bigr|\bigr) \bigr) \Biggl( \bigl|
X^{\nu}_{\cdot} - \widehat{X}{}^{\pi,\nu
}_{\cdot} \bigr| +
h + \Biggl| \int_0^T \nu_s \,ds - \sum
_{k=1}^n h \nu^{\pi
,\nu}_k
\Biggr| \Biggr).
\end{eqnarray*}
With similar arguments as used at the beginning of the proof of
Proposition \ref{prop:Y_integ}, we have
\[
\sup_{\nu\in\Uc} \E \bigl[ e^{2C (|X^{\nu}_{\cdot}| +
|\widehat
{X}{}^{\pi,\nu}_{\cdot}| )} \bigr] < +\infty\qquad
\mbox{for every } C > 0.
\]
Finally, it follows by \eqref{eq:Vd}, \eqref{eq:error_discretSemiMart}
together with the Cauchy--Schwarz inequality that~\eqref{eq:error_discretV}
holds true.
\end{pf}

\subsubsection{Invariance principle in approximation of semimartingales}

Let $X^\pi$ be a semimartingale on the discrete time grid $(t_k)_{0
\le k \le n}$ in the probability space $(\Om, \Fc, \P)$. We have also
characteristics $B^\pi$ and $M^\pi$ defined by decomposition with
respect to the natural filtration of $X^\pi$, let $\Delta A^\pi$,
$\Delta B^\pi$ be the conditional increment terms as $\Delta A^{\pi
,\nu
}$, $\Delta B^{\pi,\nu}$ defined at the beginning of Section~\ref{subsubsec:disc_semimart}. Suppose in addition that there is a constant
$C_0 > 0$ such that $\E|\Delta X^\pi_k|^3 \le C_0$ and $(\frac{1}{h}
\Delta A^\pi_k, \frac{1}{h} \Delta B^\pi_k) \in E_h \subset S_d^+
\times
\R
^d$ a.s., where $E_h$ is defined in \eqref{eq:defEh}.

Let $H \dvtx E_h \times[0,1] \to\R^d$ be a measurable mapping such that for
every $(a,b) \in E_h$ and random variable $U$ with uniform distribution
on $[0,1]$,
%
\begin{eqnarray}
\label{eq:propH}
 \E H(a,b,U) &=& b h, \qquad\operatorname{Var} H(a,b,U) = a h\quad \mbox{and}
 \nonumber
 \\[-8pt]
 \\[-8pt]
 \nonumber
  \E
\bigl|H(a,b,U)\bigr|^3 &<& C_H
\end{eqnarray}
for a constant $C_H \ge C_0$.

Now, on another probability space $(\Omb, \Fcb, \Pb)$ equipped with
$U_1, \ldots, U_n$ and $\Ub_1, \ldots, \Ub_n$ which are i.i.d. with
uniform distribution on $[0,1]$, we can approximate the distribution of
$X$ in $\Om$ by sums of random variables of the form $H(a,b,U)$ in
$\Omb$.

\begin{Lemma} \label{lemma:Xd_Xh}
There is a constant $C$ such that for every $\Theta> 0$, we can
construct two semimartingales $\Xb^\pi$ and $\Xb^h$ on $(\Omb, \Fcb
, \Pb
)$ as well as $\Delta\Ab^\pi$ and $\Delta\Bb^\pi$ such that $(\Xb
^\pi
, \Delta\Ab^\pi, \Delta\Bb^\pi)$ in $(\Omb, \Fcb, \Pb)$ has
the same
distribution as that of $(X^\pi, \Delta A^\pi, \Delta B^\pi)$ in
$(\Om,
\Fc, \P)$. Moreover,
%
\begin{eqnarray}
\label{eq:Xd_Xh}  \Xb^h_k& =& \Xb^\pi_0
+ \sum_{i=1}^k H\bigl(\Delta
\Ab^\pi_i, \Delta \Bb ^\pi_i,
\Ub_i\bigr)\quad\mbox{and}
\nonumber
\\[-8pt]
\\[-8pt]
\nonumber
\Pb \Bigl( \max_{1 \le k \le n} \bigl| \Xb^\pi_k
- \Xb^h_k \bigr| > \Theta \Bigr) &\le &C \frac{C_H n}{ \Theta^3}.
\end{eqnarray}
\end{Lemma}
We refer to Lemma 3.2 of Dolinsky \cite{Dolinsky} for a technical
proof, where the main idea is to use the techniques of invariance
principle of Sakhanenko \cite{Sakhanenko}.

\begin{Remark}
The process $\Xb^h$ is defined in \eqref{eq:Xd_Xh} by characteristics
of $\Xb^\pi$, with $\Ub$ as well as function $H$.
If we stay in the general stochastic control problem context, $H$ is a
function depending on $\mathbf{x}$,
it follows then that the process $\Xb^h$ is constructed by using
functions of the form $H(t_k, \Xb^\pi, \ldots)$,
which may not be an admissible controlled semimartingale defined in
\eqref{eq:defXh}.
This is the main reason for which we need to suppose that $H$ is
independent of $\mathbf{x}$ in \eqref{eq:condUplus} to deduce the convergence
rate in Theorem \ref{theo:convergence_rate}.
\end{Remark}
%

\begin{Remark}\label{rem:CH_value}
With $H_h$ given by \eqref{eq:defHh}, set
%
\begin{equation}
\label{eq:HequalFh} H(a,b,x):= H_h\bigl(0,0,a+b b^{\top} h, b,
x\bigr)\qquad \forall(a,b,x) \in E_h \times\R^d.
\end{equation}
Since $(a+ b b^{\top} h, b) \in E$ for every $(a,b) \in E_h$, then
function \eqref{eq:HequalFh} above is well defined.
Moreover, it follows by Lemma \ref{lemma:proba_density} that $H$
satisfies \eqref{eq:propH} with $C_H \le C_0 n^{-3/2}$, for $C_0$
independent of $n$.
In particular, let $\Theta= n^{- 1/8}$, and then
\[
\Pb \biggl( \max_{1 \le k \le n} \bigl| \Xb_k^\pi-
\Xb_k^h \bigr| > \frac
{1}{n^8} \biggr) \le
C_1 \frac{1}{n^8}
\]
for another constant $C_1$ independent of $n$ (or equivalently of $h:= T/n$).
\end{Remark}

\subsubsection{Proof of Theorem \texorpdfstring{\protect\ref{theo:convergence_rate}}{2.8}}

By Theorem \ref{theo:proba_interp} and Lemma \ref
{lemma:error_discretSemiMart}, we only need to prove separately that
%
\begin{eqnarray}
\label{eq:firstIneqRate} V^\pi&\le& V_0^h +
C_{\eps} h^{1/8- \eps},
\\
\label{eq:secondIneqRate} V_0^h &\le& V^\pi +
C_{\eps} h^{1/8- \eps}
\end{eqnarray}
and
%
\begin{equation}
\label{eq:thirdIneqRate} \bigl| V_0^h - V^\pi \bigr| \le C
h^{1/8} \qquad\mbox{if } L \mbox{ and } \Phi \mbox{ are bounded.}
\end{equation}

\paragraph*{First inequality \protect\eqref{eq:firstIneqRate}}
For every $\nu\in\Uc$ as well as the discrete-time semimartingale
$X^{\pi,\nu}$ defined at the beginning of Section~\ref{subsubsec:disc_semimart}, we can construct, following Lemma \ref
{lemma:Xd_Xh}, $(\Xb^\pi, \Delta\Ab^\pi, \Delta\Bb^\pi)$ and
$\Xb^h$
in a probability space $(\Omb, \Fcb, \Pb)$ with $H(a,b,x):=
H_h(0,0,a+b b^{\top} h, b, x)$ as in Remark \ref{rem:CH_value}, such
that the law of $(\Xb^\pi, \Delta\Ab^\pi, \Delta\Bb^\pi)$ under
$\Pb$
is the same as $(X^{\pi,\nu}, \Delta A^{\pi,\nu}, \Delta B^{\pi
,\nu})$
in $(\Om, \Fc, \P)$, and~\eqref{eq:Xd_Xh} holds true for every
$\Theta
> 0$. Fix $\Theta:= h^{1/8}$, and denote
\[
\mathcal{E}:= \Bigl\{ \max_{0\le k \le n}\bigl |\Xb^\pi_k
- \Xb ^h_k\bigr| > \Theta \Bigr\}.
\]
Denote also
\[
\nub^\pi = \bigl(\nub^\pi_k
\bigr)_{1 \le k \le n}\qquad \mbox{with } \nub ^\pi_k:=
\frac{1}{h} \bigl( \Delta\Ab^\pi_k, \Delta
\Bb^\pi_k \bigr).
\]
By the same arguments we used to prove claim \eqref
{eq:claim_exp_bound}, we know that\break  $\E^{\Pb}  [ \exp
(C(|\widehat
{\Xb}{}^\pi_{\cdot}| + |\widehat{\Xb}{}^h_{\cdot}|)  )  ] $ is
bounded by a constant independent of $\nu\in\Uc$.
It follows by the definition of $\Psi_d$ in \eqref{eq:defPsi_h} as
well as \eqref{eq:condPhi1} and \eqref{eq:condPhi2} that for every
$\eps> 0$, there is a constant $C_{\eps}$ independent of $\nu\in\Uc$
such that
%
\begin{eqnarray}
\label{eq:Rate_interm1}
&&\E^{\Pb} \bigl[\bigl | \Psi_d\bigl(
\Xb^\pi, \nub^\pi\bigr) - \Psi_d\bigl(\Xb
^h, \nub ^\pi\bigr) \bigr| \bigr]\nonumber \\
&&\qquad\le C \E^{\Pb}
\bigl[ \exp \bigl(C\bigl(\bigl|\widehat{\Xb}{}^{\pi}_{\cdot}\bigr| +\bigl |
\widehat{\Xb}{}^h_{\cdot}\bigr|\bigr) \bigr)\bigl | \widehat{
\Xb}{}^\pi_{\cdot
} - \widehat{\Xb}{}^h_{\cdot}
\bigr| \bigr]
\nonumber
\\[-8pt]
\\[-8pt]
\nonumber
&&\qquad\le C_{\eps} \bigl( h^{1/8} + \Pb(\mathcal{E})^{1/(1-8\eps)}
\bigr)
\\
&&\qquad\le C_{\eps} h^{1/8 - \eps},\nonumber
\end{eqnarray}
where the second inequality follows from H\"older's inequality and
Remark~\ref{rem:CH_value}.\vadjust{\goodbreak}

Next, we claim that
%
\begin{equation}
\label{eq:claim_Rate}  \E^{\Pb} \bigl[ \Psi_d\bigl(
\Xb^h, \nub^\pi\bigr) \bigr] = \E^{\Pb} \Biggl[
\sum_{k=0}^{n-1} h L \bigl( t_k,
\widehat{\Xb}{}^h_{\cdot}, \nub^\pi _{k+1}
\bigr) + \Phi\bigl(\widehat{\Xb}{}^h_{\cdot}\bigr) \Biggr]
 \le
V_0^h.
\end{equation}
Then by the arbitrariness of $\nu\in\Uc$, it follows by the
definition of $V^\pi$ in \eqref{eq:Vd} that~\eqref{eq:firstIneqRate}
holds true. Hence we only need to prove the claim \eqref{eq:claim_Rate}.

We can use the randomness argument as in Dolinsky, Nutz and Soner \cite
{DolinskyNutzSoner} for proving their Proposition 3.5. By the
expression of $\Xb^h$ in \eqref{eq:Xd_Xh}, using regular\vspace*{1pt} conditional
probability distribution, there is another probability space $(\tilde
\Om, \tilde\Fc, \tilde\P)$ together with independent uniformly
distributed random variables $(\tilde U_k)_{1 \le k \le n}$, $(\tilde
U_k')_{1 \le k \le n}$ and measurable functions $\Pi_k \dvtx [0,1]^k
\times
[0,1]^k \to E_h$
such that with
\[
(\Delta\tilde A_k, \Delta\tilde B_k):=
\Pi_k\bigl(\tilde U_1, \ldots, \tilde U_k,
\tilde U_1', \ldots, \tilde U_k'
\bigr)
\]
and
\[
\tilde X_k:= x_0 + \sum
_{i=1}^k H(\Delta\tilde A_i, \Delta
\tilde B_i, \tilde U_i),
\]
the distribution of $(\tilde X_k, \Delta\tilde A_k, \Delta\tilde
B_k)_{1 \le k \le n}$ in $(\tilde\Om, \tilde\Fc, \tilde\P)$ equals
to $ (\Xb^h_k, \Delta\Ab^\pi_k,\break   \Delta\Bb^\pi_k)_{1 \le k \le
n}$ in
$(\Omb, \Fcb, \Pb)$. Denote for every $u = (u_1, \ldots, u_n) \in
[0,1]^{n}$,
\[
\tilde X_k^u:= x_0 + \sum
_{i=1}^k H \bigl( \Pi_i(\tilde
U_1, \ldots, \tilde U_i, u_1, \ldots,
u_i), \tilde U_i \bigr).
\]
Since $H$ is given by \eqref{eq:HequalFh}, it follows by the
definition of $E_h$ in \eqref{eq:defEh} as well as that of $V^h_0$ in
\eqref{eq:Vh0} that, with strategy $\tilde\nu^u_k:= \frac{1}{n}
\Pi
_k(\tilde U_1, \ldots, \tilde U_k, u_1, \ldots, u_k)$,
\[
\E\bigl[ \Psi_d\bigl( \tilde X^u, \tilde
\nu^u\bigr) \bigr]\le V^h_0.
\]
And hence,
\begin{eqnarray*}
\E \bigl[ \Psi_d\bigl(\Xb^h, \nub^\pi\bigr)
\bigr] &=& \E \bigl[ \Psi_d\bigl( \tilde X^{\tilde U'}, \tilde
\nu^{\tilde U'}\bigr) \bigr] \\
&=& \int_{[0,1]^{n}} \E \bigl[
\Psi_\pi\bigl( \tilde X^u, \tilde\nu^u\bigr)
\bigr] \,du \le V^h_0.
\end{eqnarray*}
Therefore, we proved the claim, which completes the proof of
inequality~\eqref{eq:firstIneqRate}.

\paragraph*{Second inequality \protect\eqref{eq:secondIneqRate}}

Let $\tilde{F}_h(a,b,x)$ denote the distribution function of the
random variable $b h + a^{1/2} W_h$, where $W_h$ is of Gaussian
distribution $N(0,h I_d)$. Denote also by $\tilde{F}_h^{-1}(a,b,x)$ the
inverse function of $x \mapsto\tilde{F}_h(a,b,x)$. Then for every
$X^{h,\phi}$ with $\phi\in\Ac_h$, we can construct $\Xb^h$ as in
\eqref{eq:Xd_Xh} with $H(a,b,x) = \tilde{F}_h^{-1}(a,b,x)$ such that
its distribution is closed to that of $X^{h,\phi}$. By the same
arguments as in the proof of \eqref{eq:firstIneqRate}, we can prove
\eqref{eq:secondIneqRate}.

\paragraph*{Third inequality \protect\eqref{eq:thirdIneqRate}}

When $\Phi$ and $L$ are both bounded, we can improve the estimations
in \eqref{eq:Rate_interm1} to
\begin{eqnarray*}
&& \E^{\Pb} \bigl| \Psi_d\bigl(\Xb^\pi,
\nub^\pi\bigr) - \Psi_d\bigl(\Xb^h,
\nub^\pi \bigr)\bigr |
\\
&&\qquad\le 2 |\Psi_d|_{\infty} \Pb(\Ec) + C \E^{\Pb}
\bigl[ \exp \bigl( C\bigl(\bigl|\widehat{\Xb}{}^\pi_{\cdot}\bigr| + \bigl|
\widehat{\Xb}{}^h_{\cdot}\bigr|\bigr) \bigr) \bigr] h^{1/8}
\le C h^{1/8}.
\end{eqnarray*}
And all the other arguments in the proof of \eqref{eq:firstIneqRate}
and \eqref{eq:secondIneqRate} hold still true. We hence complete the
proof of \eqref{eq:thirdIneqRate}.

\section{The implementation of the scheme}
\label{sec:implem_scheme}

We shall discuss some issues for the implementation of the scheme
\eqref{eq:NumSchem}.

\subsection{The degenerate case}
\label{subsec:implem_G_expec}

The numerical scheme \eqref{eq:NumSchem} demands that $\sigma
(t,\mathbf{x}, u)
\ge\sigma_0(t,\mathbf{x}) > \eps_0 I_d$ for every $(t,\mathbf
{x},u) \in Q_T$ in
Assumption \ref{assum:G}, which implies that the volatility part should
be all nondegenerate. However, many applications are related to
degenerate cases.

\begin{Example} \label{exam:UVM}
Let $d=1$, and $E = [\underline{a}, \overline{a}]$, $\mu\equiv0$ and
$\sigma(u):= \sqrt{u}$ for $u \in E$. A concrete optimization problem
is given by
\[
\sup_{\nu\in\Uc} \E \Phi \biggl( X^{\nu}_T,
\int_0^T \nu_t \,dt \biggr).
\]
Introducing $\tilde X^{\nu}_t:= \int_0^t \nu_s \,ds$, the above problem
turns to be
\[
\sup_{\nu\in\Uc} \E \Phi \bigl( X^{\nu}_T,
\tilde X^{\nu}_T \bigr),
\]
which can be considered in the framework of \eqref{eq:V}. However, the
volatility matrix of the controlled process $(X^{\nu}, \tilde X^{\nu})$
is clear degenerate.
\end{Example}

The above example is the case of variance option pricing problem in
uncertain volatility model in finance.

\begin{Example} \label{exam:OV}
An typical example of variance option is the option ``call sharpe''
where the payoff function is given, with constants $S_0$ and $K$, by
\[
\Phi (X_T, V_T ):= \frac{  ( S_0 \exp(X_T - V_T/2) - K
)^+ }{\sqrt{V_T}}.
\]
\end{Example}

To make numerical scheme \eqref{eq:NumSchem} implementable in the
degenerate case, we can perturb the volatility matrix. Concretely,
given an optimization problem \eqref{eq:V} with coefficients $\mu$ and
$\sigma$, we set
%
\begin{eqnarray}
\label{eq:sigma_eps}  \sigma^{\eps}(t,\mathbf{x},u)&:=& \bigl(\sigma
\sigma^{\top
}(t,\mathbf{x},u) + \eps ^2 I_d
\bigr)^{1/2},
\nonumber
\\[-8pt]
\\[-8pt]
\nonumber
 a^{\eps}(t,\mathbf{x},u)&:=& \sigma^{\eps
}
\bigl(\sigma^{\eps
}\bigr)^{\top}(t,\mathbf{x},u).
\end{eqnarray}
Clearly, $a^{\eps}$ is nondegenerate.
Given $\nu\in\Uc$, let $X^{\nu,\eps}$ be the solution to SDE
%
\begin{equation}
\label{eq:SDE_eps} X^{\nu,\eps}_t = x_0 + \int
_0^t \mu\bigl(s, X^{\nu,\eps}_{\cdot},
\nu _s\bigr) \,d s + \int_0^t
\sigma^{\eps}\bigl(s, X^{\nu,\eps}_{\cdot}, \nu
_s\bigr) \,dW_s.
\end{equation}
Then a new optimization problem is given by
\[
V^{\eps}:= \sup_{\nu\in\Uc} \E \biggl[ \int
_0^T L\bigl(t, X^{\nu
,\eps
}_{\cdot},
\nu_t\bigr) \,dt + \Phi\bigl(X_{\cdot}^{\nu,\eps}\bigr)
\biggr],
\]
which is no longer degenerate.
A similar idea was also illustrated in Guyon and Henry-Labord\`ere
\cite{GuyonHenryLabordere} as well as in Jakobsen \cite{Jakobsen} for
degenerate PDEs.
We notice in addition that by applying Ito's formula on the process
$|X_t^{\nu} - \tilde X^{\nu,\eps}_t |^2$,
then taking expectations and using classical method with Gronwall's
lemma, we can easily get the error estimation
\[
\E \sup_{0 \le t \le T} \bigl| X_t^{\nu} - \tilde
X^{\nu,\eps}_t \bigr|^2 \le C \eps^2.
\]
It follows that when $L$ and $\Phi$ satisfy conditions \eqref
{eq:condPhi1} and \eqref{eq:condPhi2}, we have
\[
\bigl|V - V^{\eps}\bigr| \le C \eps.
\]

\subsection{The simulation-regression method}
To make scheme \eqref{eq:NumSchem} implementable, a natural technique
is to use the simulation-regression method to estimate the conditional
expectations arising in the scheme. First, given a function basis, we
propose a projection version of scheme \eqref{eq:NumSchem}. Next,
replacing the $L^2$-projection by least-square regression with
empirical simulations of $X^0$, it follows an implementable scheme. The
error analysis of the simulation-regression method has been achieved by
Gobet, Lemor and Warin \cite{Gobet2006} in the context of BSDE
numerical schemes. In this paper, we shall just describe the
simulation-regression method for our scheme and leave the error
analysis for further study.

\subsubsection{The Markovian setting}

In practice, we usually add new variables in the optimization problem
and make the dynamic of $X^0$ [given by \eqref{eq:defX0}] Markovian.
Suppose that for $d' > 0$, there are functions $\mub_{0,k} \dvtx \R^d
\times
\R
^{d'} \to\R^d$, $\sigmab_{0,k}\dvtx \R^d \times\R^{d'} \to S_d$ and $s_k
\dvtx \R
^d \times\R^{d'} \to\R^d $ for every $1 \le k \le n$ such that
\[
\mu_0\bigl(t_k,\widehat{X}{}^0_{\cdot}
\bigr) = \mub_{0,k}\bigl(X^0_k,
S^0_k\bigr),\qquad \sigma_0
\bigl(t_k,\widehat{X}{}^0_{\cdot}\bigr) =
\sigmab_{0,k}\bigl(X^0_k, S^0_k
\bigr)
\]
and $S^0_{k+1}:=  s_{k+1}(S^0_k, X^0_{k+1})$.
Then $(X^0_k, S^0_k)_{0 \le k \le n}$ is a Markovian process from~\eqref{eq:defX0}. Suppose further that there are\vadjust{\goodbreak} functions $(\mub_k,
\sigmab_k, \Lb_k) \dvtx \R^d \times\R^{d'} \times E \to\R^d
\times S_d \times\R$ and
$\Phib\dvtx  \R^d \times\R^{d'} \to\R$ such that
\begin{eqnarray*}
\mu\bigl(t_k,\widehat{X}{}^0_{\cdot}, u\bigr)& =&
\mub_k\bigl(X^0_k, S^0_k,
u\bigr), \qquad \sigma\bigl(t_k,\widehat{X}{}^0_{\cdot},
u\bigr) = \sigmab_k\bigl(X^0_k,
S^0_k,u\bigr)\quad\mbox{and}
\\
 L\bigl(t_k, \widehat{X}{}^0_{\cdot},u
\bigr) &= &\Lb_k\bigl(X^0_k,
S^0_k, u\bigr),\qquad \Phi \bigl( \widehat{X}{}^0_{\cdot}
\bigr) = \Phib \bigl( X^0_n, S^0_n
\bigr).
\end{eqnarray*}
Then it is clear that the numerical solution $Y^h_k$ of \eqref
{eq:NumSchem} can be represented as a measurable function of $(X^0_k,
S^0_k)$, where the function $G$ in \eqref{eq:functionG} turns to be
%
\begin{eqnarray}
\label{eq:functionGb} && \Gb(t_k,x,s,\gamma, p)
\nonumber
\\
&&\qquad:= \sup_{u \in E} \biggl( \Lb_k(x,s,u) +
\frac{1}{2} \bigl( \sigmab_k \sigmab^{\top}_k(x,s,u)
- \sigma_{0,k} \sigma^{\top}_{0,k}(x,s) \bigr) \cdot
\gamma\\
&&\hspace*{176pt}\qquad\quad{}+ \mub_k(x,s,u) \cdot p \biggr).\nonumber
\end{eqnarray}

\begin{Remark}
In finance, when we consider the payoff functions of exotic options
such as Asian options and lookback options, we can usually add the
cumulative average, or cumulative maximum (minimum) to make the system
Markovian.
\end{Remark}

\subsubsection{The projection scheme}

To simplify the notation, let us just give the scheme for the case $d=
d'= 1$, although, in general, this case can be easily deduced; we also
omit the superscript $h$ for $(Y, Z, \Gamma)$.

Let $(p^Y_{k,i})_{1 \le i \le I,  0 \le k \le n-1}$ be a family of
basis functions where every $p^Y_{k,i}$ is function defined on $\R^2$
so that
\[
\Sc_k^Y:= \Biggl\{ \sum
_{i=1}^I \alpha_i
p_{k,i}^Y\bigl(X^0_k,
S^0_k\bigr), \alpha\in\R^I \Biggr\}
\]
is a convex subclass of $L^2(\Om, \Fc_T) $. A projection operator
$\Pc
^Y_k$ is defined by
%
\begin{equation}
\label{eq:defPcY} \Pc^Y_k (U):= \arg\min
_{S \in\Sc^Y_k} \E | U - S |^2 \qquad\forall U \in
L^2(\Om, \Fc_T).
\end{equation}
Similarly, with basis functions $p^Z_{k,i}$ and $p^{\Gamma}_{k,i}$, we
can define $\Sc_k^Z$, $\Sc_k^{\Gamma}$ as well as the projections
operators $\Pc_k^Z$, $\Pc_k^{\Gamma}$. Inspired by \cite
{GobetTurd}, we
propose the following two projection schemes: with the same terminal condition
\[
\hat Y_n = \Phi\bigl(X^0_T,
S^0_T\bigr)\dvtx
\]
\paragraph*{First scheme}
\[
\cases{ \displaystyle\hat{Y}_k = \Pc^Y_k
\bigl( \hat Y_{k+1} + h \Gb\bigl(t_k, X^0_k,
S^0_k,\hat\Gamma_{k}, \hat Z_k
\bigr) \bigr), \vspace*{2pt}
\cr
\displaystyle\hat Z_k = \Pc^Z_k
\biggl( \hat Y_{k+1} \bigl(\sigma_0^{\top}
\bigr)^{-1} \frac
{\Delta W_{k+1}}{h} \biggr), \vspace*{2pt}
\cr
\displaystyle\hat
\Gamma_k = \Pc^{\Gamma}_k \biggl( \hat
Y_{k+1} \bigl(\sigma _0^{\top
}\bigr)^{-1}
\frac{\Delta W_{k+1}^{\top} \Delta W_{k+1} - h I_d}{h^2} \sigma _0^{-1} \biggr). } %
\]

\paragraph*{Second scheme}
\[
\cases{ \displaystyle\hat Y_k = \Pc^Y_k
\Biggl( \hat Y_T + \sum_{i=k}^{n-1}
h \Gb\bigl(t_k, X^0_k, S^0_k,
\hat\Gamma_i, \hat Z_i\bigr) \Biggr), \vspace*{2pt}
\cr
\displaystyle\hat Z_k = \Pc^Z_k \Biggl( \Biggl[ \hat
Y_T + \sum_{i=k+1}^{n-1} h \Gb
\bigl(t_k, X^0_k, S^0_k,
\hat\Gamma_i, \hat Z_i\bigr) \Biggr] \bigl(
\sigma_0^{\top
}\bigr)^{-1} \frac{\Delta W_{k+1}}{h}
\Biggr), \vspace*{2pt}
\cr
\displaystyle\hat\Gamma_k = \Pc^{\Gamma}_k
\Biggl( \Biggl[ \hat Y_T + \sum_{i=k+1}^{n-1}
h \Gb\bigl(t_k, X^0_k, S^0_k,
\hat\Gamma_i, \hat Z_i\bigr) \Biggr] \bigl(
\sigma_0^{\top}\bigr)^{-1}\vspace*{2pt}\cr
\hspace*{113pt}\displaystyle{}\times \frac{\Delta W_{k+1}^{\top} \Delta W_{k+1} - h
I_d}{h^2}
\sigma_0^{-1} \Biggr). } %
\]

We note that the numerical solutions are of the form $\hat Y_k =
y_k(X^0_k, S^0_k)$, $\hat Z_k = z_k(X^0_k, S^0_k)$ and $\hat\Gamma_k =
\gamma_k(X^0_k, S^0_k)$ with functions $y_k,  z_k,  \gamma_k$.
\subsubsection{Empirical regression scheme}

The simulation-regression scheme consists of simulating $M$ empirical
processes by \eqref{eq:defX0}, denoted by $(X^{0,m},\break   S^{0,m})_{1 \le m
\le M}$, and then replacing the projection of \eqref{eq:defPcY} by
empirical least square method to estimate functions $y_k$, $z_k$ and
$\gamma_k$. Concretely, with the simulation-regression method, the
first scheme turns to be
%
\begin{equation}
\label{eq:yzgamma} y_k = \sum_{i=1}^I
\widehat\alpha{}^y_{k,i} p^Y_{k,i},\qquad
z_k = \sum_{i=1}^I
\widehat\alpha{}^z_{k,i} p^Z_{k,i},\qquad
\gamma_k = \sum_{i=1}^I
\widehat\alpha{}^{\gamma}_{k,i} p^{\Gamma}_{k,i},
\end{equation}
where
\begin{eqnarray*}
\widehat\alpha{}^y_k &:=& \arg\min_{\alpha\in\R^I}
\sum_{m=1}^M \Biggl( \sum
_{i=1}^I \alpha_i
p^Y_{k,i}\bigl(X^{0,m}_k,
S^{0,m}_k\bigr)
\\
&&\hspace*{61pt}{} - y_{k+1}\bigl(X^{0,m}_{k+1},
S^{0,m}_{k+1}\bigr) - h G(t_k, \cdot, \gamma, z)
\bigl(X^{0,m}_k, S^{0,m}_k\bigr)
\Biggr)^2,
\end{eqnarray*}
and $\widehat\alpha{}^y$, $\widehat\alpha{}^{\gamma}$ are also given by
the corresponding least square method. Similarly, we can easily get an
empirical regression scheme for the second projection scheme.

Finally, we finish by remarking that in error analysis as well as in
practice, we usually need to use truncation method in formula \eqref
{eq:yzgamma} with the a priori estimations of $(\hat{Y}_k, \hat Z_k,
\hat\Gamma_k)$.

\begin{appendix}\label{app}
\section*{Appendix}
We shall give here the proofs of Lemmas \ref{lemma:weakCvg} and \ref
{lemma:equiv_WS}. The arguments are mainly due to Section~8 of Kushner
\cite{Kushner1990}, we adapt his idea of proving his Theorems 8.1, 8.2
and 8.3 in our context.

We first recall that given $\phi_h \in\Ac_h$, $X^{h,\phi_h}$ is
defined by \eqref{eq:defXh} and\break  $m^{h,\phi_h}(t, du):= \delta_{(\phi
_h)_k}(du)$ for $t \in(t_k,t_{k+1}]$. Denote $\phi^h_s:= (\phi_h)_k$
for $s \in(t_k,t_{k+1}]$ and $\E_k^{h,\phi_h}[\cdot]:= \E[\cdot|
\Fc
_k^{h,\phi_h}]$ for $\Fc_k^{h,\phi_h}:= \sigma(X^{h,\phi_h}_0,
\ldots,
X^{h,\phi_h}_k)$. Then for every $\varphi\in C_b^{\infty}(\R^d)$, it
follows by Taylor expansion that
\setcounter{equation}{0}
\begin{eqnarray}
\label{eq:Mh_martingale} && \E^{h,\phi_h}_k \bigl[ \varphi
\bigl(X^{h,\phi_h}_{k+1}\bigr) \bigr] - \varphi
\bigl(X^{h,\phi_h}_k\bigr)
\nonumber
\\[-8pt]
\\[-8pt]
\nonumber
&&\qquad= \E^{h,\phi_h}_k \biggl[ \int_{t_k}^{t_{k+1}}
\Lc^{s, \widehat{X}{}^{h,\phi_h}, \phi^h_s} \varphi\bigl(\widehat{X}{}^{h,\phi_h}_s\bigr)
\,d s \biggr] + \eps_h,
\end{eqnarray}
where $|\eps_h| \le C ( h^{3/2} +  h \rho(h))$, and $\rho$ is the
continuity module of $\mu$ and $\sigma$ in $t$ given by \eqref
{eq:cond_mu_sigma}, with a constant $C$ depending on $\varphi$ but
independent of $(h,\phi_h)$.

\begin{pf*}{Proof of Lemma \ref{lemma:weakCvg}}
(i) First, let
$(\Pb_h)_{h \le h_0}$ be the sequence of probability measures on $\Omb
^d$ given in the lemma. Suppose that $\Pb_h$ is induced by $(\widehat
{X}{}^{h,\phi_h}, m^{h,\phi_h})$ with $\phi_h \in\Ac_h$, then by
\eqref
{eq:propR}, and it is clear that there is a constant $C_3$ such that
for all $0 \le s \le t \le T$,
\[
\sup_{0 < h \le h_0} \E^{\Pb_{h}} \bigl[ \bigl|\omega^d_t
- \omega ^d_s\bigr|^3 \bigr] = \sup
_{0 < h \le h_0} \E \bigl[ \bigl| \widehat{X}_t{}^{h,\phi
_h} -
\widehat{X}_s{}^{h,\phi_h}\bigr |^3 \bigr]
 \le
C_3 |t-s|^{3/2},
\]
and hence it follows the precompactness of $ (\Pb_h|_{\Om^d}
)_{h \le h_0}$.
Further, since $E$ is supposed to be a compact Polish space, it
follows by Prokhorov's theorem that $\Mbf([0,T] \times E)$ (the space of
all positive measures $m$ on $[0,T] \times E$ such that $m([0,T]
\times E) =
T$) is compact under the weak convergence topology.
Then $\Mbf$ is also compact as a closed subset of $\Mbf([0,T] \times E)$.
It follows that the class of probability measures $(\Pb_h |_{\Mbf})_{h
\le h_0}$ on $\Mbf$ is precompact (still by Prokhorov's theorem).
Therefore, $(\Pb_h)_{ h \le h_0}$ is precompact.
Suppose that $\Pb$ is a limit measure of $(\Pb_h)_{h \le h_0}$, we
shall show that $\Pb\in\Pcb_R$. It is enough to show that for every
$\varphi\in C_b^{\infty}(\R^d)$,
%
\begin{eqnarray}
\label{eq:Mphi_martingale} \E^{\Pb_h} \bigl[ f \bigl(\omega^d_{s_i},
m_{s_i}(\psi_j), i\le I, j\le J \bigr) \bigl(
M_t(\varphi) - M_s(\varphi) \bigr) \bigr] \to0
\nonumber
\\[-8pt]
\\[-8pt]
\eqntext{\mbox
{as } h \to0}
\end{eqnarray}
for arbitrary $I$, $J \in\N$, $ s_i < s < t$, $\psi_j \in C_b([0,T]
\times E)$ and bounded continuous function $f$, where the process
$M(\varphi
)$ is defined by \eqref{eq:defM_phi}. Since $\Pb_h$ is induced by
$(\widehat{X}{}^{h,\phi_h}, m^{h,\phi_h})$ with $\phi_h \in\Ac_h$, then
\begin{eqnarray*}
&& \E^{\Pb_h} \bigl[ f \bigl(\omega^d_{s_i},
m_{s_i}(\psi_j), i\le I, j\le J \bigr) \bigl(
M_t(\varphi) - M_s(\varphi) \bigr) \bigr]
\\
&&\qquad= \E \biggl[ f \biggl( \widehat{X}{}^{h,\phi_h}_{s_i}, \int
_0^{s_i} \psi _j\bigl(
\phi^h_r\bigr) \,dr, i\le I, j\le J \biggr)\\
&&\hspace*{11pt}\qquad\quad{}\times \biggl[
\varphi\bigl(\widehat {X}{}^{h,\phi
_h}_t\bigr) - \varphi\bigl(
\widehat{X}{}^{h,\phi_h}_s\bigr)  -\int_s^t \Lc^{s, \widehat{X}{}^{h,\phi_h}, \phi^h_s} \varphi\bigl(
\widehat {X}{}^{h,\phi
_h}_r\bigr) \,dr \biggr] \biggr],
\end{eqnarray*}
which turns to $0$ as $h \to0$ by taking conditional expectations and
using~\eqref{eq:Mh_martingale}.

(ii) Suppose that $\Pb_{\delta} \in\Pcb_{S_0}$ is
induced by a
controlled process $X^{\nu^{\delta}}$ and the control $\nu^{\delta}
\in
\Uc_0$ is of the form $\nu^{\delta}_s = w(s,X_{r_i^k}^{\nu^{\delta}})$,
where $w(s,\mathbf{x}):= w_k(\mathbf{x}_{r_i^k},  i\le I_k)$ when $s
\in(\delta k,
\delta(k+1)]$ for functions $(w_k)_{k \ge0}$ and a constant $\delta> 0$.
Let $\P_{\delta}$ be the probability measure on $\Om^d$ induced by
$X^{\nu^{\delta}}$, which is clearly the unique probability measure on
$\Om^d$ under which $X_0 = x_0$ a.s. and
%
\begin{equation}
\label{eq:MartPbO} \varphi(X_t) - \int_0^t
\Lc^{s, X_{\cdot}, w(s,X_{\cdot})} \varphi (X_s) \,ds
\end{equation}
is a $\F^d$-martingale for every $\varphi\in C_b^{\infty}(\R^d)$,
where $X$ is the canonical process of $\Om^d$, and $\Lc$ is defined by
\eqref{eq:defLc}.

Now, for every $h \le h_0$, let us consider the strategy $\phi\in\Ac
_h$ defined by
\[
\phi_k(x_0, \ldots, x_k):=
w(t_k, \widehat{x}).
\]
Denote by $\P_h$ the probability measure induced by $X^{h,\phi}$ on
$\Om^d$, it follows by the same arguments as in proving \eqref
{eq:Mphi_martingale}, together with the uniqueness of solution to the
martingale problem associated to \eqref{eq:MartPbO}, that $\P_h \to
\P
_{\delta}$. Moreover, since under $\P_{\delta}$, $\mathbf{x}\mapsto
w(s, \mathbf{x}
)$ is continuous, it follows that $\Pb_h \to\Pb_{\delta}$, where
$\Pb
_h$ denotes the probability measure on $\Omb$ induced by $(\widehat
{X}{}^{h,\phi}, m^{h,\phi})$.
\end{pf*}

In preparation of the proof for Lemma \ref{lemma:equiv_WS}, we shall
introduce another subset of $\Mbf(\Omb^d)$. Let $\delta> 0$, we
consider the strategy $\nu^{\delta}_t = v_k(X_s, s \le k \delta)$ for
$t \in(\delta k, \delta(k+1)]$, where $v_k$ are measurable functions
defined on $C([0, \delta k], \R^d)$. Denote by $\Pcb_{S_c}$ the
collection of all probability measures induced by $(\widehat{X}{}^{\nu
^{\delta}},\break  m^{\nu^{\delta}})$ as in \eqref{eq:defPb}, with $\nu
^{\delta}$ of this form. Then it is clear that
\[
\Pcb_{S_0} \subset \Pcb_{S_c} \subset \Pcb_S
\subset \Pcb_R.
\]

\begin{pf*}{Proof of Lemma \ref{lemma:equiv_WS}} First, by almost
the same arguments as in Section~4 of El Karoui, H{\.u}{\.u} Nguyen and
Jeanblanc \cite{ElKaroui1987} (especially that of Theorem 4.10) that
for every $\Pb\in\Pcb_R$, there is a sequence of probability measures
$\Pb_n$ in $\Pcb_{S_c}$ such that $\Pb_n \to\Pb$, where the main idea
is using Fleming's chattering method to approximate a measure on $[0,T]
\times E$ by piecewise constant processes. We just remark that the uniform
continuity of $\mu$ and $\sigma$ w.r.t. $u$ in \eqref{eq:cond_mu_sigma}
is needed here, and the ``weak uniqueness'' assumption in their paper
is guaranteed by Lipschitz conditions on $\mu$ and $\sigma$. Then we
conclude by the fact that we can approximate a measurable function
$v_k(\mathbf{x})$ defined on $C([0, \delta k], \R^d)$ by functions
$w_k(\mathbf{x}
_{t_i}, i \le I)$ which is continuous (we notice that in Theorem 7.1 of
Kushner \cite{Kushner1990}, the author propose to approximate a
measurable function $v_k$ by functions $w_k$ which are constant on
rectangles).
\end{pf*}
\end{appendix}

\section*{Acknowledgements}
We are grateful to Nizar Touzi, J. Fr\'ed\'eric Bonnans, Nicole El
Karoui and two anonymous referees for helpful comments and suggestions.




\printaddresses


\begin{thebibliography}{28}


\bibitem{BarlesSouganidis}
\begin{barticle}[mr]
\bauthor{\bsnm{Barles},~\bfnm{G.}\binits{G.}} \AND
\bauthor{\bsnm{Souganidis},~\bfnm{P.~E.}\binits{P.~E.}}
(\byear{1991}).
\btitle{Convergence of approximation schemes for fully nonlinear second order equations}.
\bjournal{Asymptot. Anal.}
\bvolume{4}
\bpages{271--283}.
\bid{issn={0921-7134}, mr={1115933}}
\end{barticle}
\bptok{imsref}%
\endbibitem

\bibitem{BertsekasShreve}
\begin{bbook}[mr]
\bauthor{\bsnm{Bertsekas},~\bfnm{Dimitri~P.}\binits{D.~P.}} \AND
\bauthor{\bsnm{Shreve},~\bfnm{Steven~E.}\binits{S.~E.}}
(\byear{1978}).
\btitle{Stochastic Optimal Control: The Discrete Time Case}.
\bseries{Mathematics in Science and Engineering}
\bvolume{139}.
\bpublisher{Academic Press},
\blocation{New York}.
\bid{mr={0511544}}
\end{bbook}
\bptok{imsref}%
\endbibitem

\bibitem{Bonnans2004}
\begin{barticle}[mr]
\bauthor{\bsnm{Bonnans},~\bfnm{J.~Fr{\'e}d{\'e}ric}\binits{J.~F.}},
\bauthor{\bsnm{Ottenwaelter},~\bfnm{{\'E}lisabeth}\binits{{\'E}.}} \AND
\bauthor{\bsnm{Zidani},~\bfnm{Housnaa}\binits{H.}}
(\byear{2004}).
\btitle{A fast algorithm for the two dimensional {HJB} equation of stochastic control}.
\bjournal{M2AN Math. Model. Numer. Anal.}
\bvolume{38}
\bpages{723--735}.
\bid{doi={10.1051/m2an:2004034}, issn={0764-583X}, mr={2087732}}
\end{barticle}
\bptok{imsref}%
\endbibitem

\bibitem{Bouchard2004}
\begin{barticle}[mr]
\bauthor{\bsnm{Bouchard},~\bfnm{Bruno}\binits{B.}} \AND
\bauthor{\bsnm{Touzi},~\bfnm{Nizar}\binits{N.}}
(\byear{2004}).
\btitle{Discrete-time approximation and {M}onte-{C}arlo simulation of backward stochastic differential equations}.
\bjournal{Stochastic Process. Appl.}
\bvolume{111}
\bpages{175--206}.
\bid{doi={10.1016/j.spa.2004.01.001}, issn={0304-4149}, mr={2056536}}
\end{barticle}
\bptok{imsref}%
\endbibitem

\bibitem{Cheridito2007}
\begin{barticle}[mr]
\bauthor{\bsnm{Cheridito},~\bfnm{Patrick}\binits{P.}},
\bauthor{\bsnm{Soner},~\bfnm{H.~Mete}\binits{H.~M.}},
\bauthor{\bsnm{Touzi},~\bfnm{Nizar}\binits{N.}} \AND
\bauthor{\bsnm{Victoir},~\bfnm{Nicolas}\binits{N.}}
(\byear{2007}).
\btitle{Second-order backward stochastic differential equations and fully nonlinear parabolic {PDE}s}.
\bjournal{Comm. Pure Appl. Math.}
\bvolume{60}
\bpages{1081--1110}.
\bid{doi={10.1002/cpa.20168}, issn={0010-3640}, mr={2319056}}
\end{barticle}
\bptok{imsref}%
\endbibitem

\bibitem{DebrabantJakobsen}
\begin{barticle}[mr]
\bauthor{\bsnm{Debrabant},~\bfnm{Kristian}\binits{K.}} \AND
\bauthor{\bsnm{Jakobsen},~\bfnm{Espen~R.}\binits{E.~R.}}
(\byear{2013}).
\btitle{Semi-{L}agrangian schemes for linear and fully non-linear diffusion equations}.
\bjournal{Math. Comp.}
\bvolume{82}
\bpages{1433--1462}.
\bid{doi={10.1090/S0025-5718-2012-02632-9}, issn={0025-5718}, mr={3042570}}
\end{barticle}
\bptok{imsref}%
\endbibitem

\bibitem{Dolinsky}
\begin{barticle}[mr]
\bauthor{\bsnm{Dolinsky},~\bfnm{Yan}\binits{Y.}}
(\byear{2012}).
\btitle{Numerical schemes for {$G$}-expectations}.
\bjournal{Electron. J. Probab.}
\bvolume{17}
\bpages{1--15}.
\bid{doi={10.1214/EJP.v17-2284}, issn={1083-6489}, mr={2994846}}
\end{barticle}
\bptok{imsref}%
\endbibitem

\bibitem{DolinskyNutzSoner}
\begin{barticle}[mr]
\bauthor{\bsnm{Dolinsky},~\bfnm{Yan}\binits{Y.}},
\bauthor{\bsnm{Nutz},~\bfnm{Marcel}\binits{M.}} \AND
\bauthor{\bsnm{Soner},~\bfnm{H.~Mete}\binits{H.~M.}}
(\byear{2012}).
\btitle{Weak approximation of {$G$}-expectations}.
\bjournal{Stochastic Process. Appl.}
\bvolume{122}
\bpages{664--675}.
\bid{doi={10.1016/j.spa.2011.09.009}, issn={0304-4149}, mr={2868935}}
\end{barticle}
\bptok{imsref}%
\endbibitem

\bibitem{ElKaroui1987}
\begin{barticle}[mr]
\bauthor{\bsnm{El Karoui},~\bfnm{Nicole}\binits{N.}},
\bauthor{\bsnm{H{\.u}{\.u} Nguyen},~\bfnm{Du'}\binits{D.}} \AND
\bauthor{\bsnm{Jeanblanc-Picqu{\'e}},~\bfnm{Monique}\binits{M.}}
(\byear{1987}).
\btitle{Compactification methods in the control of degenerate diffusions: Existence of an optimal control}.
\bjournal{Stochastics}
\bvolume{20}
\bpages{169--219}.
\bid{doi={10.1080/17442508708833443}, issn={0090-9491}, mr={0878312}}
\end{barticle}
\bptok{imsref}%
\endbibitem

\bibitem{ElKarouiTan}
\begin{bmisc}[auto:STB|2014/02/12|12:18:25]
\bauthor{\bsnm{El Karoui},~\bfnm{N.}\binits{N.}} \AND
\bauthor{\bsnm{Tan},~\bfnm{X.}\binits{X.}}
(\byear{2013}).
\bhowpublished{Capacities, measurable selection and dynamic programming.
Preprint. Available at \url{http://www.cmapx.polytechnique.fr/\textasciitilde tan/}.}
\end{bmisc}
\bptok{imsref}%
\endbibitem

\bibitem{FahimTouziWarin}
\begin{barticle}[mr]
\bauthor{\bsnm{Fahim},~\bfnm{Arash}\binits{A.}},
\bauthor{\bsnm{Touzi},~\bfnm{Nizar}\binits{N.}} \AND
\bauthor{\bsnm{Warin},~\bfnm{Xavier}\binits{X.}}
(\byear{2011}).
\btitle{A probabilistic numerical method for fully nonlinear parabolic {PDE}s}.
\bjournal{Ann. Appl. Probab.}
\bvolume{21}
\bpages{1322--1364}.
\bid{doi={10.1214/10-AAP723}, issn={1050-5164}, mr={2857450}}
\end{barticle}
\bptok{imsref}%
\endbibitem



\bibitem{GobetTurd}
\begin{bmisc}[auto:STB|2014/02/12|12:18:25]
\bauthor{\bsnm{Gobet},~\bfnm{E.}\binits{E.}} \AND
\bauthor{\bsnm{Turkedjiev},~\bfnm{P.}\binits{P.}}
(\byear{2013}).
\bhowpublished{Linear regression MDP scheme for discrete
\mbox{BSDEs} under general conditions.
Preprint. Available at \url{http://hal.archives-ouvertes.fr/hal-00642685}.}
\end{bmisc}
\bptok{imsref}%
\endbibitem

\bibitem{GuyonHenryLabordere}
\begin{barticle}[auto:STB|2014/02/12|12:18:25]
\bauthor{\bsnm{Guyon},~\bfnm{J.}\binits{J.}} \AND
\bauthor{\bsnm{Henry-Labord{\`e}re},~\bfnm{P.}\binits{P.}}
(\byear{2011}).
\btitle{Uncertain volatility model: A Monte-Carlo approach}.
\bjournal{J. Comput. Finance}
\bvolume{14}
\bpages{37--71}.
\end{barticle}
\bptok{imsref}%
\endbibitem



\bibitem{Jakobsen}
\begin{barticle}[mr]
\bauthor{\bsnm{Jakobsen},~\bfnm{E.~R.}\binits{E.~R.}}
(\byear{2004}).
\btitle{On error bounds for approximation schemes for non-convex degenerate elliptic equations}.
\bjournal{BIT}
\bvolume{44}
\bpages{269--285}.
\bid{doi={10.1023/B:BITN.0000039390.33444.f2}, issn={0006-3835}, mr={2093506}}
\end{barticle}
\bptok{imsref}%
\endbibitem

\bibitem{Kushner1990}
\begin{barticle}[mr]
\bauthor{\bsnm{Kushner},~\bfnm{Harold~J.}\binits{H.~J.}}
(\byear{1990}).
\btitle{Numerical methods for stochastic control problems in continuous time}.
\bjournal{SIAM J. Control Optim.}
\bvolume{28}
\bpages{999--1048}.
\bid{doi={10.1137/0328056}, issn={0363-0129}, mr={1064717}}
\end{barticle}
\bptok{imsref}%
\endbibitem

\bibitem{KushnerDupuis}
\begin{bbook}[mr]
\bauthor{\bsnm{Kushner},~\bfnm{Harold~J.}\binits{H.~J.}} \AND
\bauthor{\bsnm{Dupuis},~\bfnm{Paul~G.}\binits{P.~G.}}
(\byear{1992}).
\btitle{Numerical Methods for Stochastic Control Problems in Continuous Time}.
\bseries{Applications of Mathematics (New York)}
\bvolume{24}.
\bpublisher{Springer},
\blocation{New York}.
\bid{mr={1217486}}
\end{bbook}
\bptok{imsref}%
\endbibitem

\bibitem{Gobet2006}
\begin{barticle}[mr]
\bauthor{\bsnm{Lemor},~\bfnm{Jean-Philippe}\binits{J.-P.}},
\bauthor{\bsnm{Gobet},~\bfnm{Emmanuel}\binits{E.}} \AND
\bauthor{\bsnm{Warin},~\bfnm{Xavier}\binits{X.}}
(\byear{2006}).
\btitle{Rate of convergence of an empirical regression method for solving generalized backward stochastic differential equations}.
\bjournal{Bernoulli}
\bvolume{12}
\bpages{889--916}.
\bid{doi={10.3150/bj/1161614951}, issn={1350-7265}, mr={2265667}}
\end{barticle}
\bptok{imsref}%
\endbibitem

\bibitem{PengGexpec}
\begin{bincollection}[mr]
\bauthor{\bsnm{Peng},~\bfnm{Shige}\binits{S.}}
(\byear{2007}).
\btitle{{$G$}-expectation, {$G$}-{B}rownian motion and related stochastic calculus of {I}t\^o type}.
In \bbooktitle{Stochastic Analysis and Applications}.
\bseries{Abel Symp.}
\bvolume{2}
\bpages{541--567}.
\bpublisher{Springer},
\blocation{Berlin}.
\bid{doi={10.1007/978-3-540-70847-6_25}, mr={2397805}}
\end{bincollection}
\bptok{imsref}%
\endbibitem

\bibitem{Sakhanenko}
\begin{bmisc}[auto:STB|2014/02/12|12:18:25]
\bauthor{\bsnm{Sakhanenko},~\bfnm{A.~I.}\binits{A.~I.}}
(\byear{2000}).
\bhowpublished{A new way to obtain estimates in the invariance
principle. In \textit{High Dimensional Probability II} (E. Gine, D. M. Mason and J. A. Wellner, eds.) \textit{Progr. Probab.} \textbf{47} 221--243. Birkh\"auser, Boston.}
\bid{mr={1857325}}
\end{bmisc}
\bptok{imsref}%
\endbibitem


\bibitem{STZquasisure}
\begin{barticle}[mr]
\bauthor{\bsnm{Soner},~\bfnm{H.~Mete}\binits{H.~M.}},
\bauthor{\bsnm{Touzi},~\bfnm{Nizar}\binits{N.}} \AND
\bauthor{\bsnm{Zhang},~\bfnm{Jianfeng}\binits{J.}}
(\byear{2011}).
\btitle{Quasi-sure stochastic analysis through aggregation}.
\bjournal{Electron. J. Probab.}
\bvolume{16}
\bpages{1844--1879}.
\bid{doi={10.1214/EJP.v16-950}, issn={1083-6489}, mr={2842089}}
\end{barticle}
\bptok{imsref}%
\endbibitem

\bibitem{SonerTouziZhang1}
\begin{barticle}[mr]
\bauthor{\bsnm{Soner},~\bfnm{H.~Mete}\binits{H.~M.}},
\bauthor{\bsnm{Touzi},~\bfnm{Nizar}\binits{N.}} \AND
\bauthor{\bsnm{Zhang},~\bfnm{Jianfeng}\binits{J.}}
(\byear{2012}).
\btitle{Wellposedness of second order backward {SDE}s}.
\bjournal{Probab. Theory Related Fields}
\bvolume{153}
\bpages{149--190}.
\bid{doi={10.1007/s00440-011-0342-y}, issn={0178-8051}, mr={2925572}}
\end{barticle}
\bptok{imsref}%
\endbibitem

\bibitem{Stroock1979}
\begin{bbook}[mr]
\bauthor{\bsnm{Stroock},~\bfnm{Daniel~W.}\binits{D.~W.}} \AND
\bauthor{\bsnm{Varadhan},~\bfnm{S.~R.~Srinivasa}\binits{S.~R.~S.}}
(\byear{1979}).
\btitle{Multidimensional Diffusion Processes}.
\bseries{Grundlehren der Mathematischen Wissenschaften}
\bvolume{233}.
\bpublisher{Springer},
\blocation{Berlin}.
\bid{mr={0532498}}
\end{bbook}
\bptok{imsref}%
\endbibitem

\bibitem{TanSplitting}
\begin{barticle}[mr]
\bauthor{\bsnm{Tan},~\bfnm{Xiaolu}\binits{X.}}
(\byear{2013}).
\btitle{A splitting method for fully nonlinear degenerate parabolic {PDE}s}.
\bjournal{Electron. J. Probab.}
\bvolume{18}
\bpages{1--24}.
\bid{doi={10.1214/EJP.v18-1967}, issn={1083-6489}, mr={3035743}}
\end{barticle}
\bptok{imsref}%
\endbibitem

\bibitem{Valadier}
\begin{barticle}[mr]
\bauthor{\bsnm{Valadier},~\bfnm{Michel}\binits{M.}}
(\byear{1994}).
\btitle{A course on {Y}oung measures}.
\bjournal{Rend. Istit. Mat. Univ. Trieste}
\bvolume{26}
\bpages{349--394}.
\bid{issn={0049-4704}, mr={1408956}}
\end{barticle}
\bptok{imsref}%
\endbibitem

\bibitem{Young}
\begin{bbook}[mr]
\bauthor{\bsnm{Young},~\bfnm{L.~C.}\binits{L.~C.}}
(\byear{1969}).
\btitle{Lectures on the Calculus of Variations and Optimal Control Theory}.
\bpublisher{W.~B. Saunders},
\blocation{Philadelphia}.
\bid{mr={0259704}}
\end{bbook}
\bptok{imsref}%
\endbibitem

\bibitem{Zhang2004}
\begin{barticle}[mr]
\bauthor{\bsnm{Zhang},~\bfnm{Jianfeng}\binits{J.}}
(\byear{2004}).
\btitle{A numerical scheme for {BSDE}s}.
\bjournal{Ann. Appl. Probab.}
\bvolume{14}
\bpages{459--488}.
\bid{doi={10.1214/aoap/1075828058}, issn={1050-5164}, mr={2023027}}
\end{barticle}
\bptok{imsref}%
\endbibitem

\end{thebibliography}
\end{document}